\def\DATE{\relax}
\magnification=1100
\baselineskip=14truept
\voffset=.75in
\hoffset=1truein
\hsize=4.5truein
\newdimen\hsizeGlobal
\hsizeGlobal=4.5truein%
\vsize=7.75truein
\parindent=.166666in
\pretolerance=500 \tolerance=1000 \brokenpenalty=5000

\footline={\vbox{\hsize=\hsizeGlobal\hfill{\rm\the\pageno}\hfill\llap{\sevenrm\DATE}}\hss}

\def\note#1{%
  \hfuzz=50pt%
  \vadjust{%
    \setbox1=\vtop{%
      \hsize 3cm\parindent=0pt\eightpoints\baselineskip=9pt%
      \rightskip=4mm plus 4mm\raggedright#1\hss%
      }%
    \hbox{\kern-4cm\smash{\box1}\hss\par}%
    }%
  \hfuzz=0pt
  }
\def\note#1{\relax}

\def\anote#1#2#3{\smash{\kern#1in{\raise#2in\hbox{#3}}}%
  \nointerlineskip}     

\newcount\equanumber
\equanumber=0
\newcount\sectionnumber
\sectionnumber=0
\newcount\subsectionnumber
\subsectionnumber=0
\newcount\snumber  
\snumber=0

\def\section#1{%
  \subsectionnumber=0%
  \snumber=0%
  \equanumber=0%
  \advance\sectionnumber by 1%
  \noindent{\bf \the\sectionnumber .~#1.~}%
}%
\def\subsection#1{%
  \advance\subsectionnumber by 1%
  \snumber=0%
  \equanumber=0%
  \noindent{\bf \the\sectionnumber .\the\subsectionnumber .~#1.~}%
}%
\def\prevs{\the\sectionnumber .\the\subsectionnumber .\the\snumber }

\long\def\Corollary#1{%
  \global\advance\snumber by 1%
  \bigskip
  \noindent{\bf Corollary~\prevs .}%
  \quad{\it#1}%
}%
\long\def\Lemma#1{%
  \global\advance\snumber by 1%
  \bigskip
  \noindent{\bf Lemma~\prevs .}%
  \quad{\it#1}%
}%
\def\Proof{\noindent{\bf Proof.~}}
\long\def\Proposition#1{%
  \advance\snumber by 1%
  \bigskip
  \noindent{\bf Proposition~\prevs .}%
  \quad{\it#1}%
}%
\long\def\Remark#1{%
  \bigskip
  \noindent{\bf Remark.~}#1%
}%
\long\def\Theorem#1{%
  \advance\snumber by 1%
  \bigskip
  \noindent{\bf Theorem~\prevs .}%
  \quad{\it#1}%
}%
\long\def\Statement#1{%
  \advance\snumber by 1%
  \bigskip
  \noindent{\bf Statement~\prevs .}%
  \quad{\it#1}%
}%
\def\ifundefined#1{\expandafter\ifx\csname#1\endcsname\relax}
\def\labeldef#1{\global\expandafter\edef\csname#1\endcsname{\prevs}}
\def\labelref#1{\expandafter\csname#1\endcsname}
\def\label#1{\ifundefined{#1}\labeldef{#1}\note{$<$#1$>$}\else\labelref{#1}\fi}

\def\preveq{(\the\sectionnumber .\the\subsectionnumber .\the\equanumber)}
\def\neq{\global\advance\equanumber by 1\eqno{\preveq}}

\def\ifundefined#1{\expandafter\ifx\csname#1\endcsname\relax}

\def\equadef#1{\global\advance\equanumber by 1%
  \global\expandafter\edef\csname#1\endcsname{\preveq}%
  \preveq}

\def\equaref#1{\expandafter\csname#1\endcsname}

\def\equa#1{%
  \ifundefined{#1}%
    \equadef{#1}%
  \else\equaref{#1}\fi}

\font\eightrm=cmr8%
\font\sixrm=cmr6%

\font\eightsl=cmsl8%

\font\eightbf=cmb8%

\font\eighti=cmmi8%
\font\sixi=cmmi6%

\font\eightsy=cmsy8%
\font\sixsy=cmsy6%

\font\eightex=cmex8%
\font\sixex=cmex6%
\font\fiveex=cmex5%

\font\eightit=cmti8%

\font\eighttt=cmtt8%

\font\tenbb=msbm10%
\font\eightbb=msbm8%
\font\sevenbb=msbm7%
\font\sixbb=msbm6%
\font\fivebb=msbm5%
\newfam\bbfam  \textfont\bbfam=\tenbb  \scriptfont\bbfam=\sevenbb  \scriptscriptfont\bbfam=\fivebb%

\font\tenbbm=bbm10

\font\tencmssi=cmssi10%
\font\sevencmssi=cmssi7%
\font\fivecmssi=cmssi5%
\newfam\ssfam  \textfont\ssfam=\tencmssi  \scriptfont\ssfam=\sevencmssi  \scriptscriptfont\ssfam=\fivecmssi%
\def\ssi{\fam\ssfam\tencmssi}%

\font\tenfrak=cmfrak10%
\font\eightfrak=cmfrak8%
\font\sevenfrak=cmfrak7%
\font\sixfrak=cmfrak6%
\font\fivefrak=cmfrak5%
\newfam\frakfam  \textfont\frakfam=\tenfrak  \scriptfont\frakfam=\sevenfrak  \scriptscriptfont\frakfam=\fivefrak%
\def\frak{\fam\frakfam\tenfrak}%

\font\tenmsam=msam10%
\font\eightmsam=msam8%
\font\sevenmsam=msam7%
\font\sixmsam=msam6%
\font\fivemsam=msam5%

\def\bb{\fam\bbfam\tenbb}%

\def\hexdigit#1{\ifnum#1<10 \number#1\else%
  \ifnum#1=10 A\else\ifnum#1=11 B\else\ifnum#1=12 C\else%
  \ifnum#1=13 D\else\ifnum#1=14 E\else\ifnum#1=15 F\fi%
  \fi\fi\fi\fi\fi\fi}
\newfam\msamfam  \textfont\msamfam=\tenmsam  \scriptfont\msamfam=\sevenmsam  \scriptscriptfont\msamfam=\fivemsam%
\def\msam{\msamfam\tenmsam}%
\mathchardef\leq"3\hexdigit\msamfam 36%
\mathchardef\geq"3\hexdigit\msamfam 3E%

\font\tentt=cmtt11%
\font\seventt=cmtt9%
\textfont\ttfam=\tentt
\scriptfont7=\seventt%
\def\tt{\fam\ttfam\tentt}%

\def\eightpoints{%
\def\rm{\fam0\eightrm}%
\textfont0=\eightrm   \scriptfont0=\sixrm   \scriptscriptfont0=\fiverm%
\textfont1=\eighti    \scriptfont1=\sixi    \scriptscriptfont1=\fivei%
\textfont2=\eightsy   \scriptfont2=\sixsy   \scriptscriptfont2=\fivesy%
\textfont3=\eightex   \scriptfont3=\sixex   \scriptscriptfont3=\fiveex%
\textfont\itfam=\eightit  \def\it{\fam\itfam\eightit}%
\textfont\slfam=\eightsl  \def\sl{\fam\slfam\eightsl}%
\textfont\ttfam=\eighttt  \def\tt{\fam\ttfam\eighttt}%
\textfont\bffam=\eightbf  \def\bf{\fam\bffam\eightbf}%

\textfont\frakfam=\eightfrak  \scriptfont\frakfam=\sixfrak \scriptscriptfont\frakfam=\fivefrak  \def\frak{\fam\frakfam\eightfrak}%
\textfont\bbfam=\eightbb      \scriptfont\bbfam=\sixbb     \scriptscriptfont\bbfam=\fivebb      \def\bb{\fam\bbfam\eightbb}%
\textfont\msamfam=\eightmsam  \scriptfont\msamfam=\sixmsam \scriptscriptfont\msamfam=\fivemsam  \def\msam{\msamfam\eightmsam}

\rm%
}

\def\poorBold#1{\setbox1=\hbox{#1}\wd1=0pt\copy1\hskip.25pt\box1\hskip .25pt#1}

\mathchardef\lsim"3\hexdigit\msamfam 2E%
\mathchardef\gsim"3\hexdigit\msamfam 26%

\def\bcap{\,{\textstyle\bigcap}\,}

\def\cvd{\buildrel {\rm d}\over\longrightarrow}
\def\d{\,{\rm d}}

\def\ds{\displaystyle}

\long\def\DoNotPrint#1{\relax}
\def\E{{\rm E}}
\def\eqd{\buildrel {\rm d}\over =}

\def\existsnl{\exists n\in(\Lambda,\Lambda^{1+\epsilon})\,:\,}
\def\finetune#1{#1}
\def\fixedref#1{#1\note{fixedref$\{$#1$\}$}}

\def\fn{F^\leftarrow\Bigl(1-{1\over n}\Bigr)}
\def\fnn{\fn\Bigl({N\over n}\Bigr)^{1/\alpha}}

\def\g#1{g_{[0,#1)}}
\def\gFLevy{Barbe and McCormick (2010)}
\def\Id{{\rm Id}}
\def\kl{k(\Lambda)}
\def\kvi{\Bigl(1-{N\over n}V_i\Bigr)_+^{\gamma-1}}

\def\lima{\lim_{a\to 0}}

\def\liminft{\liminf_{t\to\infty}}

\def\liml{\lim_{\Lambda\to\infty}}
\def\limn{\lim_{n\to\infty}}
\def\limsupL{\limsup_{\Lambda\to\infty}}
\def\limsupn{\limsup_{n\to\infty}}
\def\limsupt{\limsup_{t\to\infty}}
\def\limt{\lim_{t\to\infty}}
\def\limT{\lim_{T\to\infty}}

\def\LTwo{{\rm L}^2}
\def\ml{m_\Lambda}

\def\onebn{\One_{(b_n,\infty)}}

\def\oF{\overline F{}}
\def\oG{\overline G{}}

\def\prob#1{\hbox{\rm P}\{\,#1\,\}}
\def\proB#1{\hbox{\rm P}\bigl\{\,#1\,\bigr\}}
\def\Prob#1{\hbox{\rm P}\Bigl\{\,#1\,\Bigr\}}
\def\qed{{\vrule height .9ex width .8ex depth -.1ex}}
\def\ROne{{\cal R}_{1,n,N}}
\def\RTwo{{\cal R}_{2,n,N}}
\def\RThree{{\cal R}_{3,n,N}}

\def\sumROne{\sum_{i\in \ROne}}

\def\Var{{\rm Var}}

\def\boc{\note{{\bf BoC}\hskip-11pt\setbox1=\hbox{$\Bigg\downarrow$}%
         \dp1=0pt\ht1=0pt\ht1=0pt\leavevmode\raise -20pt\box1}}
\def\eoc{\note{{\bf EoC}\hskip-11pt\setbox1=\hbox{$\Bigg\uparrow$}%
         \dp1=0pt\ht1=0pt\ht1=0pt\leavevmode\raise 20pt\box1}}

\def\One{\hbox{\tenbbm 1}}

\def\calB{{\cal B}}

\def\calR{{\cal R}}

\def\MM{{\bb M}\kern .4pt}
\def\NN{{\bb N}\kern .5pt}
\def\RR{{\bb R}}

\def\UU{{\bb U}}

\def\M{{\ssi M}}


\pageno=1
\centerline{\bf RUIN PROBABILITIES IN TOUGH TIMES}
\vskip 8pt
\centerline{\bf\poorBold{$\widetilde{\hbox to 1cm{\hfill}}$}}
{\eightpoints
\centerline{\bf Part 1}
\centerline{\bf HEAVY-TRAFFIC APPROXIMATION FOR FRACTIONALLY}
\centerline{\bf INTEGRATED RANDOM WALKS IN THE DOMAIN OF ATTRACTION}
\centerline{\bf OF A NONGAUSSIAN STABLE DISTRIBUTION}
}
\bigskip

\centerline{Ph.\ Barbe$^{(1)}$ and W.P.\ McCormick$^{(2)}$}
\centerline{${}^{(1)}$CNRS {\sevenrm(UMR {\eightrm 8088})}, ${}^{(2)}$University of Georgia}

{\narrower
\baselineskip=9pt\parindent=0pt\eightpoints

\bigskip

{\bf Abstract.} Motivated by applications to insurance mathematics, we prove
some heavy-traffic limit theorems for process which encompass the fractionally
integrated random walk as well as some FARIMA processes, when the 
innovations are in the domain of attraction of a nonGaussian stable 
distribution.

\bigskip


\bigskip

\hfuzz=3pt
\noindent{\bf Keywords:} heavy traffic, ruin probability, fractional random
walk, FARIMA process, fractional L\'evy stable process.

\hfuzz=0pt
}

\bigskip\bigskip


\def\preveq{(\the\sectionnumber.\the\equanumber)}
\def\prevs{\the\sectionnumber.\the\snumber }

\section{Introduction and main result}%
The purpose of this paper is to study ruin probability when the claim process 
is nonstationary, has long range dependence, innovations in the domain
of attraction of a stable distribution and when the premiums can barely cover 
the claims; hence the title.

The motivation for such a study, beyond the development of some of the 
mathematics needed to build more realistic models for some insurance companies,
are manifold; this introduction seeks to describe them, relating the content
of this paper to various problems considered before.

To start with, the claim processes we are interested in encompass many 
nonstationnary fractional autoregressive moving average (FARIMA) processes
without prehistorical influence. 
As such, they also include the usual partial sum process. From an applied
perspective, these processes are of interest because they are part of the
standard models in time series, and their ontological justification
as aggregation of simpler processes (Granger, 1980) has some appeal in 
economics and econometrics.  From a theoretical perspective, their interest
lies in the fact that they are not Markovian, may not have stationary
solutions, may exhibit long range dependence, so may not be amenable to
the classical techniques, and yet are tractable. Consequently, a technical
understanding of these models yields a greater understanding of the underlying
stochastic phenomenon involved in simpler models. Indeed, as less technical
tools become available, we have to resort to more fundamental aspects of the
process involved. In that respect, the main contribution of this paper is 
threefold: firstly, it reveals the role of extreme values in fractional random
walks with innovations in the domain of attraction of a nonGaussian 
distribution; secondly, it gives a method of proof which, unlike all those we
are aware of for the classical random walk, does not rely on either some
form of Kolmogorov's maximal inequality or the Wiener-Hopf 
factorization ---~see a
discussion of the classical proofs in Shneer and Wachtal (2009); thirdly, 
it shows that the so-called exponential representation of uniform
order statistics may be used in sequential problems, even though this 
representation is nonsequential in nature.

A certain number of results known for the partial sum process have been
extended to its fractionally integrated version, and, more generally, to
FARIMA ones. For instance, motivated by applications in econometrics, 
Donsker's (1951) 
invariance principle, asserting the convergence of the rescaled partial sum
process to a Wiener process, has been extended to some FARIMA processes by
Philipps (1987) and Akonom and Gouri\'eroux (1987); the latter authors showed
that a fractional integral of the Wiener process, that is, a fractional
Brownian motion, may arise as limiting process ---~see Wu and Shao (2006) for
extensions and further references. In a similar spirit, Barbe and Broniatowski
(1998) extented Varadhan's (1966) large deviation result for partial sums
to some FARIMA processes ---~ see also Ghosh and Samorodnitsky (2009) for 
related results in the stationary case. In that extension, the derivative 
involved in Varadhan's action functional was replaced by a fractional 
derivative. The classical ruin estimate of Cram\'er has been partially
extended by Barbe and McCormick (2008b); and its heavy tail analogue,
Veraverbeke's (1977) Theorem 2, has been also studied in the setting of FARIMA
processes by Barbe and McCormick (2008a). The general thrust of these works
is to understand how classical results for partial sums extend to their
fractional analogue, somewhat paralleling the developments related to 
fractional Brownian motion and fractional L\'evy processes in probability.

In the classical applied probability area of queueing theory, a well studied
topic is that of the so-called heavy traffic approximations, referring to
an asymptotic analysis of the behavior of queues with traffic intensity
near one. In an insurance mathematics setting, this translates into studying a
risk process when the premium can barely cover the claims. For the partial sum
process, heavy traffic approximations are well understood and a presentation in
book form may be found in Resnick (2007) ---~see also
Shneer and Wachtel (2009) as well as Kosi\'nski, Boxma and Zwart (2010)
for further references and results, and Whitt (2002)
for examples of applications in queueing theory. One of the
purposes of this paper is to present an analogous result in the setting of
FARIMA processes.

In yet another direction which we will not pursue, heavy traffic approximation 
can be interpreted in terms of moving boundary crossing probability.
In particular, our result can give the probability that a fractional random
walk, or more generaly, a FARIMA process without prehistorical influence,
crosses a moving curved boundary.
\note{quote Doney? (Ann.\ Probab., 32, 1545; or check refs? or older work by Kesten?}

\bigskip

Throughout this paper we use the letter $c$ for a generic constant whose value
may change from one occurence to another.

We use the symbol ${}\lsim{}$ between two sequences, as in say $a_n\lsim b_n$,
to signify that $a_n\leq b_n\bigl(1+o(1)\bigr)$ as $n$ tends to infinity.

\bigskip


\section{Main result}
The processes which we will be dealing with are defined through an analytic
function $g$ on $(-1,1)$ and a distribution function $F$ on the real line.
These two pieces of data allow us to build a so-called $(g,F)$-process as
follows. Consider a sequence $(X_i)_{i\geq 1}$ of independent random variables,
all having $F$ for distribution function. We consider the series expansion of
$g$,
$$
  g(x)=\sum_{i\geq 0} g_i x^i \, .
$$
A $(g,F)$-process $(S_n)_{n\geq 0}$ is defined by $S_0=0$ and
$$
  S_n=\sum_{0\leq i<n} g_i X_{n-i} \, , \qquad n\geq 1 \, .
$$
Setting $X_i=0$ if $i$ is negative, and writing $B$ for the backward shift
operator acting on sequences, that is $BX_i=X_{i-1}$, we see that the above 
expression for $S_n$ amounts to
$$
  S_n=g(B)X_n \, .
$$

If $g(x)=1/(1-x)$, then $(S_n)$ is the partial sum process of the $X_i$.
If $g$ is a rational function continuous on $[\,-1,1\,]$, then a 
$(g,F)$-process is an ARMA one. If $g$ is $(1-\Id)^{-d}$ times a rational 
function continuous on $[\,-1,1\,]$, then a $(g,F)$-process is a FARIMA one.

Introducing the notation
$$
  \g n=\sum_{0\leq i<n} g_i \, ,
$$
we see that if $F$ has finite mean, then $S_n=g(B)(X_n-\E X_n)+\g n\E X_1$.
If the sequence $(g_n)$ is ultimately
positive and not summable, then $S_n$ drifts to $+\infty$ if $\E X_1$ is 
positive and $-\infty$
if $\E X_1$ is negative. Our heavy traffic approximation yields the limiting
behavior of $\max_{n\geq 0} S_n$ as the expected value of $X_1$ tends to
$0$ from below.
This amounts to assuming that the innovations are centered and seek the
asymptotic behavior of $\max_{n\geq 0} (S_n-a\g n)$ as $a$ tends to $0$
from above. As mentioned in the introduction, this can be interpreted as
a problem on moving boundary crossing, for the inequality $\max_{n\geq 0}
(S_n-a\g n) >x$ is equivalent to the process $S_n$ crossing the 
boundary $x+a\g n$.

A simple examination of known heavy traffic approximation results for 
the classical
random walk shows that different asymptotic behaviors are to be expected
according to the tail behavior of the innovation, and, in particular, according
to the
finiteness of the variance. The invariance principle of \gFLevy\ also
suggests that one should distinguish the cases where the sequence $(g_n)$
diverges to infinity and that where it tends to $0$, corresponding respectively
to a fractional integration and differentiation of the random walk. In this
paper, we will concentrate on the fractional integration; a companion paper 
deals with the fractional differentiation. This discussion, as well as 
technical requirements for the proof lead us to introduce some assumptions.

The summability of the $g_i$ is related to the behavior of $g$ near $-1$ and
$1$, and bears on the long range dependence properties of the process. 
We will restrict ourselves to what Granger (1988) called
the generalized integrated processes, assuming that
$$
  g(1-1/\Id) \hbox{ is regularly varying of positive index $\gamma$.}
  \eqno{\equa{assumptionGRV}}
$$

Karamata's theorem for power series asserts that if $(g_n)$ is asymptotically
equivalent to a monotone sequence, then \assumptionGRV\ is equivalent to 
$(g_n)$ being regularly varying of index $\gamma-1$. In this case, $g_i/g_n
\sim (i/n)^{\gamma-1}$ as $n$ tends to infinity and $i/n$ stays bounded away
for the origin and infinity. We will assume more; firstly, that
$$
  \hbox{$(g_n)$ is normalized regularly varying,}
  \eqno{\equa{gNRV}}
$$
meaning that (see Bingham, Goldie and Teugels, 1989, Theorem 1.9.8)
$$
  {g_{n+1}\over g_n} = 1+{\gamma-1\over n}\bigl( 1+o(1)\bigr)
$$
as $n$ tends to infinity, and, secondly, that there exists
a positive $\delta$ such that
$$
  \limn n^\delta \sup_{n^{-\delta}\leq i/n\leq 1}\Bigl|{g_i\over g_n}
  -\Bigl({i\over n}\Bigr)^{\gamma-1}\Bigr| = 0 \, .
  \eqno{\equa{hypMg}}
$$
We will prove that assumptions \gNRV\ and \hypMg\ hold for FARIMA processes;
therefore, they are not overwhelmingly restrictive in applications. Note 
that \hypMg\ implies \assumptionGRV\ and that there is no loss of generality 
to assume that $\delta<1/2$; indeed, if \hypMg\ holds for some $\delta$, then 
it holds for any smaller one. Furthermore, \gNRV\ implies \assumptionGRV\ as
well.

As far as the innovations are concerned, we assume that they have a mean but
no variance, and, more precisely, that
\setbox1=\vbox{\hsize=3in\par\noindent $F$ is centered and in the domain of 
attraction of a stable distribution of index $\alpha$ in $(1,2)$.\par}
$$
  \box1
  \eqno{\equa{assumptionStable}}
$$

Whenever $G$ is a cummulative distribution function, we write $\oG$ for
$1-G$. We write $F_*$ for the distribution function of $|X_1|$. Assumption
\assumptionStable\ implies that one of the tails of $F$ is regularly varying of
index $-\alpha$ and that $F$ is tail balanced, meaning the following. Write
$\M_{-1}F$ for the distribution function of $-X$. Then $\oF_*$ coincides with
$\oF+\overline{\M_{-1}F}$ on the positive half-line when $F$ is continuous.
The tail balance condition is that $\oF\sim p\oF_*$ and $\overline{\M_{-1}F}
\sim q \oF_*$ at infinity where $p$ and $q$ are nonnegative numbers which
add to $1$. For simplicity, we consider throughout the paper and without
mentioning it any further that $p$ does not vanish.

Writing
$$
  F^\leftarrow(u)=\inf\{\, x\,:\, F(x)\geq u\,\}
$$
for the c\`agl\`ad quantile function associated to $F$, \assumptionStable\
implies that $F^\leftarrow(1-1/\Id)$ is regularly varying of index $1/\alpha$;
if $q$ does not vanish, $(\M_{-1}F)^\leftarrow(1-1/\Id)$ is regularly varying 
of the same index $1/\alpha$. Paralleling \hypMg, we assume that for some 
positive $\kappa$,
$$
  \limt t^\kappa \sup_{t^{-\kappa}\leq\lambda\leq t^\kappa}
  \Bigl|{F^\leftarrow(1-\lambda/t)\over F^\leftarrow(1-1/t)}
  -\lambda^{-1/\alpha}\Bigr| = 0 \, .
  \eqno{\equa{hypRVRate}}
$$
If this assumption holds for some $\kappa$, then it holds for any smaller one.
While this assumption is stated in a form convenient for our usage, its
meaning is made more explicit in the following result, whose proof is deferred
to section 4.

\Proposition{\label{hypRVMeaning}
  Let $\kappa$ be a positive real number less than both $1$ and $2/(\alpha+1)$.
  The following are equivalent as $t$ tends to infinity:
  \medskip

  \noindent (i) $\ds t^\kappa \sup_{t^{-\kappa}\leq\lambda\leq t^\kappa}
  \Bigl|{\ds F^\leftarrow(1-\lambda/t)\over\ds F^\leftarrow(1-1/t)}
  -\lambda^{-1/\alpha}\Bigr| = o(1)$.

  \smallskip

  \noindent (ii) $\ds F^\leftarrow(1-1/t)=ct^{1/\alpha}
    \bigl(1+o(t^{-\kappa(1+1/\alpha)})\bigr)$.

  \medskip

  \noindent Moreover, if $F$ is continuous and increasing, 
  it is also equivalent to

  \medskip

  \noindent (iii) $\oF(t)=(c/t)^\alpha\bigl(1+o(t^{-\kappa(\alpha+1)})\bigr)$.
}

\bigskip

It is likely that our need for \hypRVRate\ is an 
artifact of the technique used in the proof, and that our result holds in a 
much greater generality. This issue is discussed after the proof,
in section \fixedref{3.6}.

Considering $(1/g_n)$, 
Proposition \hypRVMeaning\ implies that condition \hypMg\ is equivalent
to the existence of a positive $\epsilon$ such that
$$
  g_n=c n^{\gamma-1}\bigl(1+o(n^{-\epsilon})\bigr)
$$
as $n$ tends to infinity.

Concerning the lower tail we will assume either an analogue of \hypRVRate,
namely
$$
  \limt t^\kappa \sup_{t^{-\kappa}\leq\lambda\leq t^\kappa}
  \Bigl|{(\M_{-1}F)^\leftarrow(1-\lambda/t)\over (\M_{-1}F)^\leftarrow(1-1/t)}
  -\lambda^{-1/\alpha}\Bigr| = 0 \, ,
  \eqno{\equa{hypRVRateLT}}
$$
an assumption which is relevant when $q$ does not vanish and in some cases
when $q$ vanishes, or, forcing $q$ to vanish,
$$
  \overline{M_{-1}F}(t)\leq c \oF(t\log t)\log t
  \quad\hbox{ultimately.}
  \eqno{\equa{tailDominance}}
$$
Note that while assumptions \hypRVRateLT\ and \tailDominance\ do not cover
all possible distributions, it is not much more restrictive than
\hypRVRate\ in practical applications; indeed all classical distributions which
satisfy \hypRVRate\ satisfy either \hypRVRateLT\ or \tailDominance.

The distribution function $F$ yields the L\'evy measure $\nu$, defined by
its density with respect to the Lebesgue measure $\lambda$,
$$
  {\d \nu\over\d\lambda}(x)=p\alpha x^{-\alpha-1}\One_{(0,\infty)}(x)
  + q\alpha (-x)^{-\alpha-1} \One_{(-\infty,0)}(x) \, .
$$
It induces a L\'evy stable process $L_0$ with L\'evy measure $\nu$, that is
a process with selfsimilar and independent increments, such that, under
\assumptionStable,
$$
  \E e^{itL_0(1)}
  = \exp\Bigl(\int (e^{itx}-1-itx)\d\nu(x)\Bigr) \, .
$$
The subscript $0$ is to indicate that this process is centered.
A fractional L\'evy stable process is defined through the Riemann-Liouville
integral
$$
  L_0^{(\gamma-1)}(t)
  =\gamma\int_0^t (t-u)^{\gamma-1}\d L_0(u) \, .
$$
We will use the function
$$
  k={\Id\over F_*^\leftarrow(1-1/\Id)} \, .
$$
It is regularly varying of positive index $1-1/\alpha$.

Our main result is the following.

\Theorem{\label{GammaGTOne}
  Assume that $\gamma$ is greater than $1$ and that \gNRV, \hypMg, 
  \assumptionStable\ and \hypRVRate\ hold. If either \hypRVRateLT\ or 
  \tailDominance\ hold, then
  $$
    \lima {1\over ag\bigl(1-1/k^\leftarrow(1/a)\bigr)}
    \sup_{n\geq 1} (S_n-a\g n)
    \eqd \sup_{t\geq 0} \bigl(L_0^{(\gamma-1)}(t)-t^\gamma\bigr) \, .
  $$
  Moreover, the random variable $\sup_{t\geq 0} \bigl(L_0^{(\gamma-1)}(t)
  -t^\gamma\bigr)$ is almost surely finite.
}

\bigskip

\noindent{\bf Example.} We consider a FARIMA process without prehistorical 
influence defined as follows.
Let $\Theta$ and $\Phi$ be two real polynomials, the roots of $\Phi$ being
outside the complex unit disk and $\Theta$ not vanishing at $1$. The FARIMA
process
$$
  (1-B)^\gamma \Phi(B) S_n=\Theta(B)X_n
$$
is a $(g,F)$-process with $g=(1-\Id)^{-\gamma}\Theta/\Phi$. We assume that
$\gamma$ is greater than $1$. This process is then a 
fractionally integrated random walk.

Lemma 6.1 in \gFLevy\ implies that assumption \gNRV\ holds.
Checking that \hypMg\ holds is easy, because Lemmas 6.1 and 6.2 in \gFLevy\
show that when $g$ is $(1-\Id)^{-\gamma}\Theta/\Phi$, there exists a converging
sequence $(a_n)$ such that
$$
  g_n= cn^{\gamma-1} \Bigl(1+{a_n\over n}\Bigr) \, .
$$
This implies that whenever $i$ and $n$ tend to infinity with $i$ at most $n$,
$$\eqalignno{
  n^\delta\Bigl|{g_i\over g_n}-\Bigl({i\over n}\Bigr)^{\gamma-1}\Bigr|
  &{}=n^\delta \Bigl({i\over n}\Bigr)^{\gamma-1} \Bigl|{1+a_i/i\over 1+a_n/n}
    -1\Bigr| \cr
  &{}\lsim n^\delta \Bigl|{a_i\over i}-{a_n\over n}\Bigr| \, .
  &\equa{exA}\cr
  }
$$
If $i$ is in the range $[\,n^{1-\delta},n\,]$, then $n^\delta /i$ tends to
$0$ whenever $\delta$ is less than $1/2$, and, similarly, $n^\delta/n$ tends
to $0$ as $n$ tends to infinity. Since $(a_n)$ converges, we see that
\exA\ implies \hypMg.

To fix the ideas, consider a distribution function $F$ such 
that $\oF(x)\sim c x^{-\alpha}$. 
Then $F^\leftarrow(1-1/t)\sim (ct)^{1/\alpha}$ as $t$ tends to infinity.
Since $F_*^\leftarrow(1-1/t)\sim p^{-1/\alpha}F^\leftarrow(1-1/t)$, we obtain
$$
  k(t)\sim (p/c)^{1/\alpha} t^{(\alpha-1)/\alpha} 
$$
as $t$ tends to infinity. It follows that
$$
  k^\leftarrow(1/a)\sim (c/p)^{1/(\alpha-1)} a^{-\alpha/(\alpha-1)}
$$
as $a$ tends to $0$. Since $g(1-1/x)\sim x^\gamma \Theta(1)/\Phi(1)$ as
$x$ tends to infinity, we have
$$
  a g\Bigl(1-{1\over k^\leftarrow(1/a)}\Bigr)
  \sim (c/p)^{\gamma/(\alpha-1)} a^{(\alpha-1-\alpha\gamma)/(\alpha-1)} 
  {\Theta(1)\over\Phi(1)}
$$
as $a$ tends to $0$. Hence, assuming that \hypRVRate\ holds, we obtain that
$$
  \sup_{n\geq 0} (S_n-a\g n) 
  \sim  (c/p)^{\gamma/(\alpha-1)}{\Theta(1)\over\Phi(1)} 
  a^{1-\gamma\alpha/(\alpha-1)} 
  \sup_{t\geq 0} \bigl( L_0^{(\gamma-1)}(t)-t^\gamma\bigr)
$$
as $a$ tends to $0$. In particular, the left hand side grows like
$a^{1-\gamma\alpha/(\alpha-1)}$. It is interesting to note that the exponent
involved depends on $\alpha$ and $\gamma$ only 
through $\gamma\alpha/(\alpha-1)$, that is $\gamma$ times the conjugate
exponent of $\alpha$.

\bigskip


\section{Proof of Theorem \GammaGTOne}
For the classical random walk, there exists two ways of proving a heavy
traffix approximation: one based on the Wiener-Hopf factorization proposed by
Kingman (1961, 1962, 1965), one based on a functional limit theorem proposed
by Prohorov (1963). We follow Prohorov's approach.

Throughout the proof we will use many times the following form of Karamata's 
theorem for power series (see Bingham, Goldie and Teugels, 1989, Corollary 
1.7.3). If $(g_n)$ is regularly varying of positive index $\gamma-1$, it is
asymptotically equivalent to an increasing sequence and
$$
  g_n
  \sim {\gamma \g n\over n}
  \sim {g(1-1/n)\over n\Gamma(\gamma)}
  \eqno{\equa{Karamata}}
$$
as $n$ tends to infinity.

To proceed with the proof,
up to an asymptotic equivalence, define $\Lambda=\Lambda(1/a)$ by the relation
$$
  a\kl\sim 1
  \eqno{\equa{LambdaDef}}
$$
as $a$ tends to $0$. It follows from Barbe and McCormick's (2010) Theorem 5.2 
that, in the sense of weak$*$ convergence of distribution of stochastic 
processes in ${\rm D}[\,0,\infty)$ endowed with the topology of uniform 
convergence on compactas,
$$
  {\kl \over \g\Lambda}S_{\lfloor\Lambda\Id\rfloor}
  \cvd L_0^{(\gamma-1)}
$$
as $\Lambda$ tends to infinity. Since \Karamata\ and \hypMg\ imply that 
$(\g n)$ is a regularly varying sequence of index $\gamma$, \LambdaDef\ 
implies that we have the convergence of stochastic processes
$$
  {\kl\over \g \Lambda} (S_{\lfloor\Lambda\Id\rfloor}-a\g{\Lambda\Id})
  \cvd L_0^{(\gamma-1)}-\Id^\gamma \, .
  \eqno{\equa{invariance}}
$$
Consequently, for any positive $T$,
$$
  {\kl\over \g\Lambda}\sup_{0\leq n\leq\Lambda T} (S_n-a\g n)
  \cvd \sup_{0\leq t\leq T} \bigl(L_0^{(\gamma-1)}(t)-t^\gamma\bigr)
$$
as $a$ tends to $0$.

Since $\sup_{0\leq t\leq T} \bigl(L_0^{(\gamma-1)}(t)-t^\gamma\bigr)$ is
nondecreasing in $T$, it converges almost surely as $T$ tends to infinity, 
possibly to infinity. Hence, to prove Theorem \GammaGTOne, it suffices to 
show that
$$
  \limT\limsup_{a\to 0}\prob{ \exists n>\Lambda T \,:\, S_n>a\g n}
  = 0
  \eqno{\equa{ThATightA}}
$$
and $\sup_{t\geq 0} \bigl(L_0^{(\gamma-1)}(t)-t^\gamma\bigr)$ is almost surely
finite.

As pointed out by
Shneer and Wachtel (2009), or differently in Szczotka and Woykzy\'nski (2003),
the main difficulty in proving a heavy traffic approximation for sums
is to show that the maximum of the process does not occur
too far in time, that is, in our case proving \ThATightA.
In the context of $(g,F)$-processes this task is far more involved
than for ordinary random walks, mostly because there is no analogue
of Kolmogorov's maximal inequality. In order to explain our proof, that is,
the remainder of this paper, we need to make a preliminary study of \ThATightA.

Given \LambdaDef, we see that \ThATightA\ is equivalent to
$$
  \limT\limsupL\prob{\exists n>\Lambda T \,:\, S_n>\g n/\kl}
  =0 \, .
  \eqno{\equa{ThATightB}}
$$
Substituting $\Lambda$ for $\Lambda T$, using that
$$
  k(\Lambda/T)\sim T^{(1/\alpha)-1}\kl
$$
as $\Lambda$ tends to infinity, and using that $1-1/\alpha$ is positive,
substituting $T$ for $T^{1-1/\alpha}$, we obtain that \ThATightB\ is equivalent
to
$$
  \limT\limsupL \prob{\exists n>\Lambda \, :\, S_n>T\g n/ \kl}=0
  \, .
  \eqno{\equa{ThATightC}}
$$

We can now explain how to prove Theorem \GammaGTOne. The proof has four 
main steps. The first two aim at showing that instead of considering
all $n$ exceeding $\Lambda$ in \ThATightC, we can reduce the range to
all $n$ between $\Lambda$ and $\Lambda^{1+\epsilon}$, where $\epsilon$ is 
positive but can be chosen as one wishes. This is achieved by showing in 
the first step that the innovations coming from the central part of the
distribution can be ignored. In the second step, a simple bound
on the contribution of the largest innovations permits us to show that if
the event involving $S_n$ in \ThATightC\ occurs, it is very likely that
$n$ is less than $\Lambda^{1+\epsilon}$. Being able to concentrate on
the range of $n$ between $\Lambda$ and $\Lambda^{1+\epsilon}$, the third
step consists in showing a result similar to that of the first step, namely
that the innovations not too large can be ignored; while this mid-range
depends on $n$ in the first step, this dependence will be, in some sense,
less so in the third step. The fourth step, by far the most complicated, 
consists in showing that the contribution of the extreme innovations to $S_n$,
properly rescaled, can be approximated by a fractional L\'evy process, 
uniformly in the range of $n$ between $\Lambda$ and $\Lambda^{1+\epsilon}$.
While this shares some similarity with \invariance\ and will be proved with
a technique inspired by our proof of \invariance, this is more difficult
than proving \invariance. The reason is that in \invariance\ we approximate
the process $S_n$ on $[\,0,\Lambda\,]$, which is $\Lambda$ times the fixed
compact set $[\,0,1\,]$. In contrast in our problem, the 
set $[\,\Lambda,\Lambda^{1+\epsilon}\,]$ should be thought as $\Lambda$ times
$[\,1,\Lambda^\epsilon\,]$, that is $\Lambda$ times an interval whose
length diverges with $\Lambda$. This forces us to develop a sequential
analogue of the representation used in \gFLevy, and it is likely that the
technique used will be of value for related boundary crossing problems. 
A fifth step concludes the proof, mostly taking care of the lower tail
of the distribution and doing some bookeeping.

\bigskip

\def\prevs{\the\sectionnumber .\the\subsectionnumber.\the\snumber }
\def\preveq{(\the\sectionnumber .\the\subsectionnumber.\the\equanumber)}

\subsection{Step 1}
We first consider the part of the process $(S_n)$ made from the not
too large innovations. In order to set up the proper thresholdings, let
$(a_n)$ and $(b_n)$ be two sequences of real numbers such that
$$
  \limn a_n=-\infty
  \qquad\hbox{and}\qquad
  \limn b_n=+\infty \, .
$$
Since $F$ obeys a tail balance condition and we suppose that $p$ does not
vanish, it is convenient to assume that
$$
  \limn b_n/(-a_n) \quad\hbox{is positive or infinite.}
$$
We define the variance of the truncated innovations,
$$
  \sigma_n=\Var\bigl(X\One_{[a_n,b_n]}(X)\bigr) \, ,
$$
and the centered and standardized `middle' innovations,
$$
  Z_{i,n}={X_i\One_{[a_n,b_n]}(X_i)-\E X\One_{[a_n,b_n]}(X)
           \over \sigma_n} \, , \qquad 1\leq i\leq n \, .
$$
The `middle' part of $S_n$ is then
$$
  M_n=\sigma_n\sum_{0\leq i<n} g_i Z_{n-i,n} \, .
$$

Given that we use $(a_n)$ and $(b_n)$ to truncate the innovations and
that we will use a quantile transformation, it is convenient, as in
\gFLevy, to take both sequences as quantiles. Since both sequences are
constructed in similar fashion, we explain that of $(b_n)$. We
consider a sequence $(\widetilde m_n)$ which is regularly varying of
positive index $\beta$ less than $1$.  We set $\widetilde
b_n=F^\leftarrow(1-\widetilde m_n/n)$. Setting $m_n=n\oF(\widetilde
b_n)$, the sequence
$$
  \hbox{$(m_n)$ is regularly varying of index $\beta$}
$$
as well. We then take
$$
  b_n=F^\leftarrow(1-m_n/n) \, .
$$
Since $(m_n)$ is a regularly varying sequence of index $\beta$, the sequence
$(b_n)$ is regularly varying of index $(1-\beta)/\alpha$.

This construction ensures that $(m_n)$ is regularly varying of index $\beta$
and $1-m_n/n$ is in the range of $F$. This ensures that the 
inequality $F^\leftarrow(1-u)>b_n$ is equivalent to $u<m_n/n$ 
(see Shorack and Wellner, 1986, \S 1, pp.\ 5--7). It is implicit that a 
similar construction is made for $(a_n)$, switching the tails.

In order to avoid heavy subscripts and many integer parts brackets, we will 
sometimes 
use the function $m(\cdot)$ defined by $m(x)=m_{\lfloor x\rfloor}$.
We will also write $m_x$ for $m_{\lfloor x\rfloor}$.

Our next proposition asserts that the middle part can be neglected in our
problem. Recall that the parameter $\beta$ regulates the growth of our 
truncation sequence used to define $M_n$.

\Proposition{\label{middleNeglect}
  For any $\beta$ positive and less than $1$, for any positive $T$,
  $$
    \liml 
    \prob{\exists n\geq\Lambda \, : \, M_n>T\g n/\kl}
    = 0 \, . 
  $$
}

\Proof Lemma 2.1.1 in \gFLevy\ asserts that 
$$
  \sigma_n\sim c b_n\sqrt{\oF(b_n)}\sim c F^\leftarrow(1-m_n/n)\sqrt{m_n/n}
$$
as $n$ tends to infinity. Inequality (2.2.1) in \gFLevy\ implies 
that for any 
positive integer $r$ there exists a constant $c_r$ such that for any 
positive $n$,
$$\eqalign{
  \Bigl|\E\Bigl({M_n\over\sigma_n\sqrt n}\Bigr)^r\Bigr|
  &{}\leq {c_r\over n} \sum_{1\leq i\leq n} |g_i|^r \cr
  &{}\sim c_r |g_n|^r \int_0^1 u^{r(\gamma-1)}\d u \, ,\cr}
$$
the asymptotic equivalence being as $n$ tends to infinity.
Using Markov's inequality and \Karamata, for any positive integer $r$,
$$\eqalignno{
  \prob{M_n>T\g n/k(\Lambda)}
  &{}\leq cg_n^r \Bigl({\kl\sigma_n\sqrt n\over T\g n}\Bigr)^r \cr
  &{}\leq c \Bigl( {\kl\over T n} F^\leftarrow(1-m_n/n)\sqrt{m_n}\Bigr)^r 
   \cr
  &{}\sim {c\over T^r} m_n^{-r/2} \Bigl({\kl\over k(n/m_n)}\Bigr)^r
  &\equa{middleNeglectA} \cr
  }
$$
as $\Lambda$ tends to infinity and uniformly in $n\geq \Lambda$. Since $\alpha$
is less than $2$, let $\eta$ be a positive real number so that $(1/2)-
(1/\alpha)+\eta$ is negative. Potter's bound implies
$$
  k\Bigl({n\over m_n}\Bigr)
  \gsim m_n^{(1/\alpha)-1-\eta}k(n)\, .
$$
It follows that \middleNeglectA\ is of order at most
$$
  {c\over T^r} m_n^{r((1/2)-(1/\alpha)+\eta)} \Bigl({\kl\over k(n)}\Bigr)^r
  \, .
$$
For $T$ and $\Lambda$ large enough, this bound is less than $1$ since 
$(1/2)-(1/\alpha)+\eta$ is negative, $r$ is positive, and $k$ is regularly 
varying of positive index. Applying Bonferroni's inequality, we obtain
$$\displaylines{\quad
  \prob{\exists n\geq\Lambda \, :\, M_n>T\g n/\kl}
  \hfill\cr\noalign{\vskip 3pt}\hfill
  \leq {c\over T} \ml^{r((1/2)-(1/\alpha)+\eta)} \kl^r
  \sum_{n\geq\Lambda} k(n)^{-r} \, .
  \quad\equa{middleNeglectB}
  \cr}
$$
Taking $r$ greater than $\alpha/(\alpha-1)$ ensures that the series of generic
term $k(n)^{-r}$ converges and that $\sum_{n\geq\Lambda} k(n)^{-r}$ is of order
$\Lambda k(\Lambda)^{-r}$
as $\Lambda$ tends to infinity. In this case, \middleNeglectB\ is of order
$\ml^{r((1/2)-(1/\alpha)+\eta)}\Lambda$. This bound is regularly varying in
$\Lambda$ of index $\beta r\bigl((1/2)-(1/\alpha)+\eta\bigr)+1$. This index is
negative whenever $r$ is large enough.\hfill\qed

\bigskip


\subsection{Step 2} We consider the contribution of the extreme innovations
to $S_n$,
$$
  T_n^+=\sum_{0\leq i<n} g_i X_{n-i}\onebn(X_{n-i}) \, .
$$

In order to understand precisely the role of \hypRVRate\ we introduce a slight
refinment. Let $\xi$ be an ultimately increasing slowly varying function,
diverging to infinity at infinity,
such that $\xi(n^2)\sim c\xi(n)$ as $n$ tends to infinity, and
$$
  \sum_{n\geq 1} {1\over n\xi(n)}<\infty  \, .
  \eqno{\equa{xiDefA}}
$$
One could take $\xi(x)$ to be $(\log x)^{1+\eta}$ or 
$(\log x)(\log\log x)^{1+\eta}$ for some positive $\eta$; 
one may simply replace $\xi(n)$ by $\log^2 n$ when reading the 
remainder of the proof. However, having introduced the function $\xi$ will 
allow us to understand the role of \hypRVRate. Sometimes we will write
$\xi_n$ instead of $\xi(n)$. The key requirement, \xiDefA,
is equivalent to the assertion that the smallest of $n$ independent random
variables uniformly distributed over $[\,0,1\,]$ is greater than $1/(n\xi_n)$
almost surely for $n$ large enough (see Geffroy, 1958, 1959, and Kiefer, 
1972); the other conditions are of technical nature. 

We introduce the following variant of \hypRVRate: there exists a real 
number $\rho$ greater than $1$ such that
$$
  \limsupn m_n^\rho 
  \sup_{1/m_n\leq\lambda\leq m_n} \Bigl|{F^\leftarrow(1-\lambda/n)
  \over F^\leftarrow(1-1/n)}-\lambda^{-1/\alpha}\Bigr| <\infty \, ,
  \eqno{\equa{hypRVRateAlt}}
$$
as well as the condition
$$
  \liminf_{n\to\infty} m_n/\xi_n>0 \, .
  \eqno{\equa{hypMnUpsilon}}
$$
Note that if $(m_n)$ is regularly varying of positive index $\beta$ less than
$\kappa$, then \hypRVRate\ implies \hypRVRateAlt\ and \hypMnUpsilon; 
indeed, \hypRVRate\ the implies
$$
  \limn n^\kappa \sup_{1/m_n\leq\lambda\leq m_n}
  \Bigl| {F^\leftarrow(1-\lambda/n)\over F^\leftarrow(1-1/n)}-
          \lambda^{-1/\alpha}\Bigr| = 0 \, ,
$$
and, considering the indices of regular variation, we can take $\rho$ to be any
number greater than $1$ and less than $\kappa/\beta$. In what
follows will rely solely on the combination \hypRVRateAlt\ and \hypMnUpsilon,
and, except specified otherwise, not on the positivity of the index of
regular variation $\beta$ of $(m_n)$. In particular, if $\beta$ is allowed to
vanish, $(m_n)$ is allowed to be slowly varying.
Ultimately, this will inform us on the role of \hypRVRate. A further discussion is
presented in section \fixedref{3.6}.

\Proposition{\label{stepTwo}
  Let $\epsilon$ be a positive real number. If
  $$
    \beta<\Bigl({\epsilon\over 1+\epsilon}{\alpha-1\over\alpha}\Bigr) 
        \wedge \kappa\, .
  $$
  then for any positive $T$,
  $$
    \lim_{\Lambda\to\infty}
    \prob{ \exists n\geq \Lambda^{1+\epsilon}\,:\, |T_n^+-\E T_n^+|
    \geq T\g n/\kl} = 0 \, .
  $$
}

In order to prove Proposition \stepTwo, we need the following estimates
on the expectation $\mu_n^+=\E X\One_{[b_n,\infty)}(X)$.

\Lemma{\label{mun}
  If \hypRVRateAlt\ holds, then
  $$
    \limsupn \Bigl|k(n)\mu_n^+-{\alpha\over\alpha-1}m_n^{1-1/\alpha}\Bigr|
    <\infty \, .
  $$
}

\Proof Let $U$ be a uniform random variable over $[\,0,1\,]$. Thinking of the 
random variable $X$ as $F^\leftarrow(1-U)$,
$$\eqalign{
  {n\mu_n^+\over F^\leftarrow(1-1/n)}
  &{}=n\int_0^{m_n/n}{F^\leftarrow(1-u)\over F^\leftarrow(1-1/n)}\d u \cr
  &{}=\int_0^{m_n} {F^\leftarrow(1-v/n)\over\ds F^\leftarrow(1-1/n)} 
   \d v\, . \cr
  }
$$
Assumption \hypRVRateAlt\ yields
$$\eqalign{
  \int_{1/m_n}^{m_n} {F^\leftarrow(1-v/n)\over F^\leftarrow(1-1/n)} \d v
  &{}=\int_{1/m_n}^{m_n} v^{-1/\alpha}+O(m_n^{-1}) \d v \cr
  &{}={\alpha\over\alpha-1} m_n^{1-1/\alpha} +O(1) \, . \cr
  }
$$
Furthermore, Potter's bound implies that for any $\eta$ positive less than
$(\alpha-1)/\alpha$, as $n$ tends to infinity,
$$
  \int_0^{1/m_n} {F^\leftarrow(1-v/n)\over F^\leftarrow(1-1/n)}\d v
  \leq 2\int_0^{1/m_n} v^{-(1/\alpha)-\eta} \d v = O(1)
$$
This proves the lemma.\hfill\qed

\bigskip

Note that Lemma \mun\ implies, as $n$ tends to infinity,
$$
  \mu_n^+\sim {\alpha\over\alpha-1} {m_n^{1-1/\alpha}\over k(n)}
  \eqno{\equa{munEquiv}}
$$

\bigskip

\noindent{\bf Proof of Proposition \stepTwo.} The proof has two steps. In the
first one we prove that, almost surely, $T_n^+$ cannot 
exceed $\g n/k(\Lambda)$ whenever $n$ exceeds
$\Lambda^{1+\epsilon}$ and $\Lambda$ is large enough. In the second one, we
prove a similar assertion on the expectation $\E T_n^+$. Recall that since
$(g_n)$ is regularly varying of positive index, it is asymptotically 
equivalent to a nondecreasing sequence.

\noindent{\it Step 1.} Let $(U_i)_{i\geq 1}$ be a 
sequence of independent random variables, uniformly distributed on $[\,0,1\,]$.
There is no loss of generality in assuming that $X_i=F^\leftarrow(1-U_i)$. 
Since we use the c\`agl\`ad version of the quantile function and $1-m_n/n$
is in the range of $F$, the inequality
$F^\leftarrow(1-U)>F^\leftarrow(1-m_n/n)$ occurs if and only if $U<m_n/n$
(see Shorack and Wellner, 1986, \S I.1, pp.5--7).  Therefore, writing $\UU_n$
for the empirical distribution function of $(U_i)_{1\leq i\leq n}$, we
have, for any $n$ large enough,
$$\eqalignno{
  T_n^+
  &{}=\sum_{0\leq i<n} g_i F^\leftarrow(1-U_{n-i})\One\{\, U_i\leq m_n/n\,\}\cr
  &{}\leq 2 g_n F^\leftarrow (1-U_{1,n}) n\UU_n(m_n/n) \, .
  &\equa{stepTwoA}\cr
  }
$$
From Theorem 1 in Kiefer (1972) we deduce that $U_{1,n}\geq 1/n\xi_n$ 
almost surely for $n$ large enough, and from Theorem 2 in Shorack and 
Wellner (1978), we conclude that $\UU_n\leq \xi_n\Id$ almost surely 
for $n$ large enough. Therefore, using \Karamata, \stepTwoA\ is 
ultimately at most
$$
  2 g_n F^\leftarrow\Bigl(1-{1\over n\xi_n}\Bigr) m_n\xi_n 
  \sim c \g n {m_n\xi_n^2\over k\bigl(n\xi_n)}\, .
  \eqno{\equa{stepTwoB}}
$$
Recall that $m(\cdot)$ is the function such that $m(x)=m_{\lfloor x\rfloor}$.
Provided $\beta$ is less than $1-1/\alpha$, the function 
$m_n\xi_n^2/k(n\xi_n)$ is regularly varying in $n$ of negative
index $\beta-1+(1/\alpha)$. Thus, \stepTwoB\ is at most
$$
  c \g n { m(\Lambda^{1+\epsilon}) \xi^2(\Lambda^{1+\epsilon})
           \over k\bigl(\Lambda^{1+\epsilon}\xi(\Lambda^{1+\epsilon})\bigr)}
$$
in the range of $n$ at least $\Lambda^{1+\epsilon}$ and for any $\Lambda$
large enough.
To see that this is less than $\g n/k(\Lambda)$, note that, considering the
index of regular varition, the inequality
$$
  { m(\Lambda^{1+\epsilon}) \xi^2(\Lambda^{1+\epsilon}) 
    \over k\bigl(\Lambda^{1+\epsilon}\xi(\Lambda^{1+\epsilon})\bigr)}
  \leq {1\over\kl}
$$
holds since
$$
  (1+\epsilon)\Bigl(\beta-1+{1\over\alpha}\Bigr)
  < -1+{1\over\alpha}
$$
whenever
$$
  \beta<{\epsilon\over 1+\epsilon}{\alpha-1\over\alpha} \, .
$$

\noindent{\it Step 2.} Using \munEquiv,
$$
  \E T_n^+\sim {\alpha\over\alpha-1}\g n {m_n^{1-1/\alpha}\over k(n)}
  \eqno{\equa{stepTwoBb}}
$$
as $n$ tends to infinity. Potter's bound to compare $k(n)$ and $k(n\xi_n)$
shows that \stepTwoBb\ is less than \stepTwoB, and therefore less 
than $\g n/\kl$ in our range of $n$ and $\Lambda$ of interest.\hfill\qed

\bigskip

Combining Propositions \middleNeglect\ and \stepTwo, we see that,
ignoring for the time being the lower tail of $F$, Theorem \GammaGTOne\ will 
be proved if we show that for some $\epsilon$
positive,
$$
  \limT\limsupL \prob{ \existsnl |T_n^+-\E T_n^+|\geq T\g n/\kl} = 0 \, .
  \eqno{\equa{stepTwoBB}}
$$
As can be seen in the remainder of the proof, the fact that the innovations
are kept in $T_n^+$ if they exceed a threshold $b_n$ in which $m_n$ 
depends on $n$ creates some 
complications. So, it is better to backtrack from \stepTwoBB, and 
instead, still using Propositions \middleNeglect\ and \stepTwo, to argue
that, still ignoring the problem of the lower tail of $F$ for the time being, 
Theorem \GammaGTOne\ can be proved by showing that for any
positive $\epsilon$,
$$
  \limT\limsupL \prob{ \existsnl |S_n|\geq T\g n/\kl} = 0 \, .
  \eqno{\equa{stepTwoC}}
$$
In the next step we show that we can now consider only the extreme values
larger than some $b_{n,\Lambda}$ calculated with a sequence $\ml$ instead 
of $m_n$.

\bigskip

\subsection{Step 3} We now concentrate on the range of $n$ between $\Lambda$
and $\Lambda^{1+\epsilon}$. Imitating the notation used in step 1, let
$$
  b_{n,\Lambda} = F^\leftarrow\Bigl(1-{\ml\over n}\Bigr)
$$
and similarly for $(a_{n,\Lambda})$. We set
$$
  \sigma_{n,\Lambda}
  =\Var\bigl(X\One_{[a_{n,\Lambda},b_{n,\Lambda}]}(X)\bigr) \, .
$$
Using the same notation as in step 1, but for a slightly different
quantity ---~note indeed that we substitute $a_{n,\Lambda}$ and $b_{n,\Lambda}$
for $a_n$ and $b_n$ ~--- 
consider the standardized middle innovations,
$$
  Z_{i,n}={X_i\One_{[a_{n,\Lambda},b_{n,\Lambda}]}(X_i)
           -\E X\One_{[a_{n,\Lambda},b_{n,\Lambda}]}(X)
           \over \sigma_{n,\Lambda}} \, , \qquad 1\leq i\leq n \, .
$$
Again, with a slight change of notation compared to step 1, the corresponding 
middle part of $S_n$ is then
$$
  M_n=\sigma_{n,\Lambda}\sum_{0\leq i<n} g_i Z_{n-i,n} \, .
$$

As in step 1, this middle part can be neglected in our
problem.

\Proposition{\label{middleNeglectThree}
  For any $\beta$ positive and less than $1$, for any positive $\epsilon$ 
  and $T$,
  $$
    \liml 
    \prob{\existsnl M_n>T\g n/\kl}
    = 0 \, . 
  $$
}

\Proof The same estimates as in Proposition \middleNeglect\ give the analogue
of \middleNeglectA, namely, that for any positive integer $r$,
$$\eqalignno{
  \prob{M_n>T\g n/k(\Lambda)}
  &{}\sim {c\over T^r} \ml^{-r/2} \Bigl({\kl\over k(n/\ml)}\Bigr)^r
  &\equa{middleNeglectThreeA} \cr
  }
$$
as $\Lambda$ tends to infinity and uniformly in $n$ 
in $(\Lambda,\Lambda^{1+\epsilon})$. Let $\eta$ be a positive real number.
Potter's bound implies that  uniformly in the
range $n$ in $(\Lambda,\Lambda^{1+\epsilon})$, as $\Lambda$ tends to infinity,
$$
  k\Bigl({n\over \ml}\Bigr)
  \gsim \ml^{(1/\alpha)-1-\eta}k(n)\, .
$$
If follows that \middleNeglectThreeA\ is of order at most
$$
  {c\over T^r} \ml^{r((1/2)-(1/\alpha)+\eta)}
  \, .
$$
Applying Bonferroni's inequality, we obtain
$$\displaylines{\quad
  \prob{\existsnl M_n>T\g n/\kl}
  \hfill\cr\hfill
  \leq {c\over T^r} \ml^{r((1/2)-(1/\alpha)+\eta)} \Lambda^{1+\epsilon} \, .
  \quad\equa{middleNeglectThreeB}
  \cr}
$$
This bound is regularly varying in
$\Lambda$ of index $\beta r\bigl((1/2)-(1/\alpha)+\eta\bigr)+1+\epsilon$. This 
index is negative whenever $\eta$ is small enough and $r$ is large 
enough.\hfill\qed

\bigskip


\def\onebn{\One_{(b_{n,\Lambda},\infty)}}

\subsection{Step 4} We consider the contribution of the extreme innovations
to $S_n$. We consider $b_{n,\Lambda}=F^\leftarrow(1-\ml/n)$ as in the 
previous subsection. With again a slight change of notation compared to step 2,
we set
$$
  T_n^+=\sum_{0\leq i<n} g_i X_{n-i}\onebn(X_{n-i}) \, .
$$
Paralleling what we did in step 2, we seek to prove the following
proposition.

\Proposition{\label{stepFour}
  Let $\epsilon$ be a positive real number. If $\beta$ is small enough,
  then for any positive $T$,
  $$
    \limT\limsupL\prob{ \existsnl |T_n^+-\E T_n^+|
    \geq T\g n/\kl} = 0 \, .
  $$
}

The proof of this proposition is far more difficult than that of Proposition
\stepTwo. It first requires several approximations of $T_n^+$. Our goal at the
end of these approximations amounts to be able to replace $T_n^+-\E T_n^+$ by
about $L_0^{(\gamma-1)}(n)/k(n)$.

Given a positive $\epsilon$ as in Proposition \stepFour, we introduce
for notational simplicity
$$
  N=\Lambda^{1+2\epsilon} \, .
$$
The exponent could as well be taken to be $1+\epsilon$ but adding an extra
$\epsilon$ will slightly simplify some of our arguments.

Following the construction in \gFLevy, let $X_{1,N}\leq X_{2,N}\leq\ldots\leq
X_{N,N}$ be the order statistics of the innovations $(X_i)_{1\leq i\leq N}$.
Let $\tau$ be the random permutation of $\{\, 1,2,\ldots,N\,\}$ such that
$$
  X_{\tau(i)} = X_{N-i+1,N} \, .
$$
We set $g_i$ to be $0$ if $i$ is negative. For any $n$ positive at most $N$, 
the equality
$$
  T_n^+=\sum_{1\leq i\leq N} g_{n-\tau(i)} X_{N-i+1,N} \onebn(X_{N-i+1,N})
  \One\{\, \tau(i)\leq n\,\}
$$
holds.
Let $(V_i)_{1\leq i\leq N}$ be a sequence of independent random variables
uniformly distributed over $[\,0,1\,]$, independent 
of $(X_{i,N})_{1\leq i\leq N}$. Let $G_N$ be their empirical distribution 
function,
$$
  G_N(x)=N^{-1}\sum_{1\leq i\leq N} \One\{\, V_i\leq x\,\} \, .
$$
Without any loss of generality, even if $F$ is not continuous, we assume 
that $\tau(i)=NG_N(V_i)$, giving
$$\displaylines{\qquad
  T_n^+= \sum_{1\leq i\leq N} g_{n-NG_N(V_i)} X_{N-i+1,N}
  \hfill\cr\hfill
  \onebn(X_{N-i+1,N})
  \One\{\, NG_N(V_i)\leq n\,\} \, .\cr}
$$
Let $(\omega_i)_{i\geq 1}$ be a sequence of independent random variables having
a standard exponential distribution. For any positive integer $i$ we define 
the partial sum $W_i=\omega_1+\cdots +\omega_i$. 
Since $(W_i/W_{N+1})_{1\leq i\leq N}$ has
the same distribution as the order statistics of $N$ independent uniform random
variables (see Shorack and Wellner, 1986, chapter 8, \S 2),
$$
  (X_{N-i+1,N})
  \eqd \biggl(F^\leftarrow\Bigl(1-{W_i\over W_{N+1}}\Bigr)
       \biggr)_{1\leq i\leq N} \, .
$$
Since we use the c\`agl\`ad version of the quantile function, the inequality
$F^\leftarrow(u)>b_{n,\Lambda}$ is equivalent to $u<\ml /n$. Therefore, introducing the
random set
$$
  \ROne=\Bigl\{\,i\,:\, {W_i\over W_{N+1}}
  \leq {\ml\over n}\, ;\,  G_N(V_i) \leq {n\over N}  \,\Bigr\}
$$
and the random variable
$$
  T_{1,n,N}^+ 
  = \sumROne g_{n-NG_N(V_i)} F^\leftarrow\Bigl( 1-{W_i\over W_{N+1}}\Bigr)
  \, ,
  \eqno{\equa{TOnePlus}}
$$
we have $(T_n^+)_{1\leq n\leq N}\eqd (T_{1,n,N}^+)_{1\leq n\leq N}$.

\bigskip

\Remark The main reason we proved Propositions \middleNeglect\ and \stepTwo\ 
is that once we can reduce $n$ to be bounded from above by some given quantity
---~in our case $\Lambda^{1+\epsilon}$~--- we can use the representation
of the innovations with $(W_i)$ and obtain \TOnePlus.

\bigskip

Let $\mu_{n,\Lambda}^+$ be $\E X\onebn(X)$, so that 
$\E T_n^+=\g n \mu_{n,\Lambda}^+$. Our discussion shows that in order 
to prove Proposition \stepFour\ it suffices to prove that
\finetune{\hfuzz=9pt}
$$\displaylines{
  \limT\limsup_{\Lambda\to\infty} 
  \Prob{ \exists n\in (\Lambda,\Lambda^{2+\epsilon}) \,:\,
  |T_{1,n,N}^+-\g n\mu_{n,\Lambda}^+|>T{\g n\over\kl} }
  \hfill\cr\hfill
  {}= 0 \, .
  \qquad\equa{stepThreeA}\cr}
$$
\hfuzz=0pt
This will be achieved by approximating $T_{1,n,N}^+$ by a much simpler process.
While the main approximation scheme follows that in \gFLevy, a main difference
lies in the sequential nature of the event involved in \stepThreeA; indeed,
contrary to \gFLevy\ we need more than an approximation of $T_{1,n,N}^+$
valid for a range of $n$ of comparable order, but for $n$ between the different
orders $\Lambda$ and $\Lambda^{1+\epsilon}$, that is over a much larger range.
To make this approximation, it is essential to have some understanding of the
set $\ROne$.

The point process $\sum_{i\geq 1} \delta_{(W_i,V_i)}$ is a Poisson process
of intensity the Lebesgue measure on $[\,0,\infty)\times [\,0,1\,]$. 
Viewing this point process in a $(w,v)$-plane, we can think of $\ROne$ as a
region in this plane given by 
$$
  \{\, (w,v) \, :\, w\leq W_{N+1} \ml/n\,;\, G_N(v)\leq n/N\,\}\, .
$$ 
So we will write about the `set' $\ROne$ when we think
of it as a subset of the integer and about the `region' $\ROne$ when we view
it as a set of points in the $(w,v)$-plane. The key fact to understand, which
we will formalize in our next result, is that, considering $\ROne$, if $n$
is small the inequality $G_N(V_i)\leq n/N$ will select few points and
because $V_i$ will be quite small, we will wait for a long time, that is we
will need to have $i$ large, in order to hit this $V_i$; but for large $i$,
the sum $W_i$ will be large and so $W_i/W_{N+1}$ will not be less than $\ml/n$.
If $n$ is large, then the inequality $W_i/W_{N+1}\leq \ml/n$ forces $i$ to be
quite small, and so, regardless of $n$, we should have very few points in 
$\ROne$. A different argument, less informative in our context though, is that
since we retained the innovations exceeding $F^\leftarrow(1-\ml/n)$ 
in $T_{1,n,N}^+$, there should be about $\ml$ of those contributing 
to $T_{1,n,N}^+$. Thus, the cardinality of $\ROne$ should be about $\ml$.

Still considering $\ROne$, if we replace $W_{N+1}$ by its near expected 
value $N$ and approximate $G_N$ by its limit, the identity
function on $[\,0,1\,]$, we should have
$$
  V_i\leq {n\over N} \qquad\hbox{and}\qquad W_i\leq \ml{N\over n}  \, .
$$
These two inequalities force $(V_i,W_i)$ to belong to the much simpler
region $\{\, (v,w)\,:\, v\leq \ml/w\,\}$, which is the subgraph of an 
hyperbola. The replacement of $G_N$ by its limit cannot be done at
this early stage and we will have to settle for less, bounding $G_N$ by a 
multiple of its limit. This leads us to introduce, for any positive $c$,
the region
$$
  \calR_{\Lambda,c}
  =\Bigl\{\, (v,w)\,:\, v\leq c{\ml\over w}\; ;\; w\leq 2N\,\Bigr\}
  \, .
$$
Our next result asserts that provided $c$ is large enough, it is very likely
for the regions $\ROne$ to be included in $\calR_{\Lambda,c}$, and that the 
minimum of the set $\ROne$ is very likely to be at least $N/n\xi_N$.

\Lemma{\label{boundRegion}
  For any positive $\eta$ there exists a real number $c$ such that for 
  any $\Lambda$ large enough,
  $$
    {\rm P}\Bigl(\bcap_{\Lambda\leq n\leq \Lambda^{1+\epsilon}} \{\,\ROne
    \subset \calR_{\Lambda,c}\,\}\Bigr) \geq 1-\eta \, .
  $$
  Moreover, viewing $\ROne$ as a set of integers,
  $$
    \liml {\rm P}\Bigl(\bcap_{\Lambda\leq n\leq\Lambda^{1+\epsilon}}
    \Bigl\{ \min\ROne\geq {N\over n\xi_N}\,\Bigr\}\Bigr) = 1 \, .
  $$
}

\Proof Let $\eta$ be an arbitrary positive real number. Inequality (2) in
Shorack and Wellner (1978) (see also Shorack and Wellner, 1986, chapter 10,
\S 3, inequality 1, p.\ 415) and the strong law of large numbers
ensure the existence of $c$ such that the event
$$
  \Omega_N = \bcap_{1\leq i\leq N} \Bigl\{\, V_{i,N}\leq c{i\over N}\,\Bigr\}
  \cap \{\, W_{N+1}\leq 2N\,\}
$$
has probability at least $1-\eta$ whenever $N$ is large enough. On $\Omega_N$,
if $i$ belongs to $\ROne$,
$$
  V_i
  \leq V_{n,N}
  \leq c {n\over N}
  \leq c \ml {W_{N+1}\over W_i} {1\over N}
  \leq 2 c {\ml\over W_i} \, ,
$$
the first equality coming from the condition $G_N(V_i)\leq n/N$, the second
inequality from being in $\Omega_N$, the third from the condition $W_i/W_{N+1}
\leq \ml/n$, and the last from being in $\Omega_N$. Thus, up to replacing
$c$ by twice as much, $\ROne\subset\calR_{\Lambda,c}$ on $\Omega_N$, proving 
the first assertion.

Given \xiDefA, Geffroy (1958/1959) or Kiefer's (1972) Theorem 1 
imply that for $c$ small enough the event 
$$
  \Omega=\bcap_{i\geq 1} \{\, V_{1,i}\geq c/i\xi_i \,\}
$$
has probability at least $1-\eta$. If $i$ belongs to $\ROne$ and $\Omega_N$
occurs, $V_i\leq c n/N$, and, in particular, $V_{1,i}\leq c n/N$. Therefore,
on $\Omega\cap\Omega_N$ we must have $1/(i\xi_i)\leq cn/N$, that is, 
$i\xi_i\geq cN/n$. In particular, $i\geq cN/(n\xi_i)$. However, if $i$ is in
$\ROne$, then $i$ is at most $N$, and since $(\xi_n)$ is ultimately monotone, 
we obtain that $i$ is at least $cN/\bigl(n\xi_N\bigr)$.\hfill\qed

\bigskip

As a consequence of Lemma \boundRegion, we can show that the region $\ROne$ 
cannot contain too many points if $\beta$ is small.

\Lemma{\label{cardinalityRegion}
  $\max_{\Lambda\leq n\leq\Lambda^{2+\epsilon}} \sharp\ROne = O_P(\ml \log N)$.
  }

\bigskip

\Proof Since $V_i$ is uniform 
over $[\,0,1\,]$ the region $\calR_{\Lambda,c}$ can be restricted to
$$
  \{\, (v,w)\,:\, v\leq (c\ml/w)\wedge 1 \, ;\, w\leq 2N\,\} \, .
$$
The Lebesgue measure of this region is of order $c \ml\log (2N)$.
The result follows since $\sum_{i\geq 1}\delta_{(V_i,W_i)}$ is
a homogenous Poisson process with mean intensity $1$ and Lemma \boundRegion\
holds.\hfill\qed

\bigskip

We will use the following lemma, which we state now for convenience.

\Lemma{\label{Poisson}
  Let $(\Pi_{i,n})_{1\leq i\leq n}$ be some Poisson random variables, possibly
  dependent, having respective means $(\lambda_{i,n})_{1\leq i\leq n}$, such
  that for some $\epsilon$ positive, 
  $\max_{1\leq i\leq n} \lambda_{i,n}=o(n^{-\epsilon})$. Then, for $k\geq
  2/\epsilon$,
  $$
    \limn\prob{\max_{1\leq i\leq n} \Pi_{i,n}\geq k}=0 \, .
  $$
}

\Proof Chernoff's inequality yields for any positive $k$,
$$
  \prob{\Pi_{i,n}\geq k}
  \leq \exp(-k\log k+k\log\lambda_{i,n}+k-\lambda_{i,n}) \, .
$$
Given the assumption on $(\lambda_{i,n})$ this upper bound is, for any $n$ 
large enough, at most $\exp(-k\log k-\epsilon k\log n+k)$. The result
follows from Bonferroni's inequality.\hfill\qed

\bigskip

Since we will make repeated use of the following simple argument or obvious
variants of it, we state it as a lemma.

\Lemma{\label{trick} 
  Let $(\epsilon_n)$ be a bounded sequence of positive real numbers.
  There exists a positive $T$ such that for any $n$ at least $\Lambda$ and
  any $\Lambda$ large enough,
  $$
    F^\leftarrow(1-1/n)g_n \epsilon_n\leq T\g n/k(\Lambda) \, .
  $$
}

\Proof Given \Karamata, the inequality amounts 
to $\epsilon_n<c k(n)/k(\Lambda)$. Since the function $k$ is regularly 
varying of positive index, 
$$
  \lim_{\Lambda\to\infty}\inf_{n\geq\Lambda} k(n)/k(\Lambda)=1 \, ,
$$ 
and the result follows.\hfill\qed

\bigskip

Having made these observations on $\ROne$, we can start a long string of
approximations of $T_{1,n,N}^+$. Referring to \TOnePlus, we first replace
$F^\leftarrow(1-W_i/W_{N+1})$ by 
$F^\leftarrow(1-1/n)(nW_i/W_{N+1})^{-1/\alpha}$. Define
$$
  T_{2,n,N}^+
  = F^\leftarrow\Bigl(1-{1\over n}\Bigr) 
  \Bigl({W_{N+1}\over n}\Bigr)^{1/\alpha} \sum_{i\in\ROne} g_{n-NG_N(V_i)}
  W_i^{-1/\alpha} \, .
$$
Our next lemma shows that we can replace $T_{1,n,N}^+$ by $T_{2,n,N}^+$ in 
order to prove \stepThreeA. Recall that $\beta$ refers to the index of
regular variation of $(\ml)$ as a function of $\Lambda$, and that, 
except if specified otherwise, we allow it
to vanish under \hypRVRateAlt\ and \hypMnUpsilon.

\Lemma{\label{approxTwo}
  If $\beta$ is less than $\kappa$, then
  $$
    \limT\limsupL\Prob{\existsnl |T_{1,n,N}^+-T_{2,n,N}^+|>T{\g n\over\kl}}
    =0\, .
  $$
}

\Proof Let $i$ be an integer in $\ROne$. We have $nW_i/W_{N+1}\leq \ml$. 
Lemma \boundRegion\ implies that except on a set whose probability can be 
made arbitrary small by taking $N$ large enough, $i$ is at 
least $N/n\xi_N$. Hence, the strong law of large numbers yields 
that $W_i$ is at least $N/2n\xi_N$ provided $N$ is large enough. Then,
$nW_i/W_{N+1}$ is at least $1/4\xi_N$. Consequently, \hypRVRateAlt\ 
and \hypMnUpsilon\ yield
$$\displaylines{\quad
  F^\leftarrow\Bigl(1-{W_i\over W_{N+1}}\Bigr)
  \hfill\cr\hfill
  {}= F^\leftarrow\Bigl(1-{1\over n}\Bigr) 
  \Bigl({nW_i\over W_{N+1}}\Bigr)^{-1/\alpha} 
  +F^\leftarrow\Bigl(1-{1\over n}\Bigr)  O_P(\ml^{-\rho}) \, ,\quad\cr
  }
$$
the $O_P(\ml^{-\rho})$ being uniform in $i$ in $\ROne$. Since $(g_n)$ is 
asymptotically equivalent to a nondecreasing sequence and
Lemma \cardinalityRegion\ holds, we obtain
$$\eqalign{
  T_{1,n,N}^+
  &{}=T_{2,n,N}^+ +F^\leftarrow\Bigl(1-{1\over n}\Bigr)g_n \ml^{-\rho}
    \sharp\ROne O_P(1)\cr
  &{}= T_{2,n,N}^+ + F^\leftarrow\Bigl(1-{1\over n}\Bigr)g_n o_P(1)\, , 
   \cr}
$$
the $o_P$-term being uniform in $n$ between $\Lambda$
and $\Lambda^{1+\epsilon}$. A variation on Lemma \trick\ implies
the result.\hfill\qed

\bigskip

In $T_{2,n,N}^+$, we replace $W_{N+1}$ by $N$, setting
$$
  T_{3,n,N}^+
  =\fnn\sum_{i\in\ROne} g_{n-NG_N(V_i)} W_i^{-1/\alpha} \, .
$$

\Lemma{\label{approxThree}
  If $\beta$ is less than $1/2$,
  $$
    \limT\limsupL\Prob{\existsnl |T_{2,n,N}^+-T_{3,n,N}^+| > T{\g n\over\kl}}
    =0 \, .
  $$
}

\Proof Taylor's formula, the central limit theorem and the strong law of large
numbers yield
$$
  W_{N+1}^{1/\alpha}-N^{1/\alpha}
  = N^{(1/\alpha)-(1/2)} O_P(1)
$$
as $N$ tends to infinity. Lemma \boundRegion\ and the strong law of large
numbers imply that if $i$ is in $\ROne$, then $W_i^{-1/\alpha}\leq 2
\bigl( n\xi_N/N\bigr)^{1/\alpha}$ whenever $N$ is large enough. 
Therefore, using Lemma \cardinalityRegion,
$$\displaylines{
  |T_{2,n,N}^+-T_{3,n,N}^+|
  \hfill\cr\hfill
  \eqalign{
  {}\leq{}& \fnn N^{-1/2} g_n\sharp \ROne 
            \Bigl({n\xi_N\over N}\Bigr)^{1/\alpha} O_P(1) \cr
  {}={}   & \fn g_n {\ml\xi^{1/\alpha}(N)\over \sqrt N} \log N O_P(1) \, , \cr}
  \cr}
$$
uniformly in $n$ between $\Lambda$ and $\Lambda^{1+\epsilon}$. The result
follows by a simple adaptation of Lemma \trick.\hfill\qed

\bigskip

Seeking to replace $\ROne$ by the slightly simpler set
$$
  \RTwo=\Bigl\{\, i\,:\, W_i\leq \ml{N\over n}\,;\, G_N(V_i)\leq {n\over N}\,
        \Bigr\} \, ,
$$
we set
$$
  T_{4,n,N}^+=F^\leftarrow\Bigl(1-{1\over n}\Bigr) 
  \Bigl({N\over n}\Bigr)^{1/\alpha} \sum_{i\in\RTwo} g_{n-NG_N(V_i)} 
  W_i^{-1/\alpha} \, .
$$

\Lemma{\label{approxFour}
  For any $\beta$ less than $(1/2)-\epsilon$,
  $$
    \limT\limsupL \Prob{\existsnl |T_{3,n,N}^+-T_{4,n,N}^+|
    \geq T{\g n\over k(\Lambda)}}
    = 0 \, .
  $$
}

\Proof If $i$ belongs to the symmetric difference $\ROne\bigtriangleup\RTwo$,
then $G_N(V_i)\leq n/N$, and either
$$
  \ml {W_{N+1}\over n}\leq W_i\leq \ml {N\over n}
  \qquad\hbox{or}\qquad
  \ml {N\over n}\leq W_i\leq \ml{W_{N+1}\over n} \, .
$$
Thus, $W_i$ lies in a region of width
$$
  {\ml\over n} |W_{N+1}-N|={\ml\over n} \sqrt N O_P(1)
$$
with an endpoint given by $\ml N/n$. In particular, $i\sim \ml N/n$ 
and $(V_i,W_i)$ lies in a region of area of order at most
$$
  {\ml\over n} \sqrt N O_P(1)
  = \ml\Lambda^{-(1/2)+\epsilon} O_P(1) \, ,
$$
the $O_P(1)$-term being uniform in $\Lambda\leq n\leq \Lambda^{1+\epsilon}$. 
In particular, taking $\beta$ less than $(1/2)-\epsilon$,
this area tends to $0$ at algebraic rate. Applying Lemma \Poisson, there exists
a positive $k$ such that
$$
  \liml \proB{\max_{\Lambda\leq n\leq \Lambda^{1+\epsilon}}\sharp
  (\ROne\bigtriangleup\RTwo)\geq k}=0 \, .
$$
Using again that all $i$ in $\ROne\bigtriangleup\RTwo$ are asymptotically
equivalent to $\ml N/n$ and using also the strong law of large numbers, we
obtain
\hfuzz=1pt
$$\eqalign{
          &\hskip -12pt|T_{3,n,N}^+-T_{4,n,N}^+|\cr
  \hskip -3pt
  {}\leq{}& \fnn g_n\hskip -1pt\max_{i\in\ROne\bigtriangleup\RTwo}\hskip -3pt W_i^{-1/\alpha}
    \sharp (\ROne\bigtriangleup\RTwo)\cr
  \hskip -3pt
  {}\leq{}& \fn g_n \ml^{-1/\alpha}O_P(1) \, . \cr
  }
$$
\hfuzz=0pt
The result then follows from Lemma \trick.\hfill\qed

\bigskip

Considering $T_{4,n,N}^+$, we seek to replace $g_{n-NG_N(V_i)}$ by 
$g_{\lfloor n-NV_i\rfloor}$. For simplicity, we will write $g_{n-NV_i}$ for
the latter. Therefore, we define
$$
  T_{5,n,N}^+=\fnn \sum_{i\in\RTwo} g_{n-NV_i}W_i^{-1/\alpha} \, .
$$

\Lemma{\label{approxFive}
  For any $\beta$ less than $\bigl((\gamma-1)\wedge 1\bigr)/4$,
  $$
    \limT\limsupL \Prob{\existsnl |T_{4,n,N}^+-T_{5,n,N}^+|>T 
    {\g n\over k(\Lambda)}} =0 \, .
  $$
}

\Proof Let $\epsilon_1$ and $\eta$ be two positive real numbers. Since the 
function $\bigl(\Id(1-\Id)\bigr)^{(1/2)-\eta}$ is a Chibisov-O'Reilly function
(see for instance Cs\"org\H o, Cs\"org\H o, Horv\`ath and Mason, 1986,
Theorem 4.2.3),
$$
  \max_{1\leq i\leq N}
  {\sqrt N |G_N(V_i)-V_i|\over \bigl(V_i(1-V_i)\bigr)^{(1/2)-\eta}}
  = O_P(1)
$$
as $N$ tends to infinity. If $i$ belongs to $\RTwo$, then $G_N(V_i)\leq n/N$
and, using Shorack and Wellner's (1978, inequality (2)) 
linear bounds, $V_i\leq c n/N$ with probability
at least $1-\epsilon_1$ provided $c$ is large enough. Thus,
$$
  \max_{i\in \RTwo} N|G_N(V_i)-V_i|\leq \sqrt N 
  \Bigl({n\over N}\Bigr)^{(1/2)-\eta} O_P(1)
  \eqno{\equa{approxFiveA}}
$$
where the $O_P(1)$-term is uniform in $n$ between $\Lambda$ 
and $\Lambda^{2+\epsilon}$. 

For any integer $r$ let
$$
  \Omega_n(r)=\max_{1\leq i\leq n-r} |g_{i+r}-g_i| \, .
$$
Inequality \approxFiveA\ implies that with probability at least $1-\epsilon_1$
provided $c$ is large enough,
$$
  \max_{i\in\RTwo}|g_{n-NG_N(V_i)}-g_{n-NV_i}|
  \leq \max_{0\leq r\leq c n^{(1/2)-\eta}N^\eta} \Omega_n(r) \, . 
  \eqno{\equa{approxFiveB}}
$$
Recall that \gNRV\ holds.
According to whether $\gamma$ is at least $2$ or not,  Lemmas 5.6 and 5.8 
in \gFLevy\ imply that whenever $\theta$ is a positive real
number less than
$(\gamma-1)\wedge 1$, the right hand side of \approxFiveB\ is at most
$cg_n(n^{-(1/2)-\eta}N^\eta)^\theta$. Since $n$ is at least $\Lambda$ in our 
range of interest, the right 
hand side of \approxFiveB\ is of order $g_n$ times $\Lambda$ at the power 
$-\theta/2+O(\eta)$. Thus, if $\eta$ is small and $n$ and $N$ are large enough,
enough,
$$\displaylines{
  |T_{4,n,N}^+-T_{5,n,N}^+|
  \hfill\cr\hfill
  \eqalign{
  {}\leq{}& \fnn g_n N^{-\theta/4} \sum_{i\in\RTwo} W_i^{-1/\alpha}O_P(1) \cr
  {}\leq{}& \fnn g_n N^{-\theta/4} \max_{i\in\RTwo}W_i^{-1/\alpha}\sharp\RTwo
            O_P(1)\, . \cr
  }
  \cr}
$$
Lemma \boundRegion\ with $\RTwo$ substituted for $\ROne$, the strong law 
of large numbers applied to the sums
$(W_i)_{i\geq 1}$ and Lemma \cardinalityRegion\ with $\RTwo$ substituted for
$\ROne$ shows that the above upper bound is, in probability, of order
$$\displaylines{\qquad
  \fnn g_n N^{-\theta/4} \Bigl({n\xi_N\over N}\Bigr)^{1/\alpha} \ml \log N
  \hfill\cr\hfill
  {}=\fn g_n N^{-\theta/4} \xi_N^{1/\alpha} \ml \log N\, .
  \qquad\cr}
$$
Thus, if $\beta$ is less than $\theta/4$, Lemma \trick\ implies the 
result.\hfill\qed

\bigskip

Next, using the regular variation of the sequence $(g_n)$, we would like 
to replace $g_{n-NV_i}$ in $T_{5,n,N}^+$ by $g_n (1-NV_i/n)^{\gamma-1}$.
This leads us to define
$$
  T_{6,n,N}^+=\fnn g_n \sum_{i\in \RTwo} 
  \Bigl( 1-{N\over n} V_i\Bigr)_+^{\gamma-1} W_i^{-1/\alpha} \, .
$$

With regard to the next lemma, recall that $\delta$ was introduced in \hypMg.

\Lemma{\label{approxSix}
  If $\beta$ is less than $\delta$, then
  $$
    \limT\limsupL \Prob{\existsnl |T_{5,n,N}^+-T_{6,n,N}^+|>T{\g n\over\kl}}
    =0 \, .
  $$
}

\Proof To approximate $T_{5,n,N}^+$ by $T_{6,n,N}^+$, we need to rely on 
assumption $\hypMg$. With respect to this assumption, we see that $\RTwo$
contains `good' points, for which
$$
  n^{-\delta}\leq {n-NV_i\over n}\leq 1\, ,
$$
guaranteeing that with \hypMg\ we can substitute $g_n(1-NV_i/n)^{\gamma-1}$
for $g_{n-NV_i}$,
and `bad' points, for which either
$$
  0\leq {n-NV_i\over n}\leq n^{-\delta} 
  \qquad\hbox{or}\qquad V_i>n/N \, .
$$
We call $\calB_{1,n,N}$ the set of all bad points for which $0\leq n-NV_i
\leq n^{1-\delta}$ and $\calB_{2,n,N}$ the set of those for which $V_i>n/N$.

Let $i$ be in $\calB_{1,n,N}$. Since it belongs to $\RTwo$, 
Lemma \boundRegion\ shows that with probability
arbitrarily close to $1$ provided $N$ is large enough, $i\geq N/n\xi_N$. 
Therefore, since $i$ is in $\RTwo$, with probability arbitrarily close to $1$
provided $N$ is large enough,
$$
  {N\over 2n\xi_N}\leq W_i\leq \ml{N\over n} \, .
$$
And since $i$ is a bad point in $\calB_{1,n,N}$,
$$
  {n-n^{1-\delta}\over N}\leq V_i\leq {n\over N}\, .
$$
When thinking of $\calB_{1,n,N}$ as a region as we did with $\ROne$, we thus
have with high probability,
$$
  \calB_{1,n,N}
  \subset \Bigl\{\, (v,w) \,:\, {n-n^{1-\delta}\over N} \leq v\leq {n\over N}
  \, ;\, {N\over 2 n\xi_N}\leq w\leq \ml {N\over n}\,\Bigr\} \, .
$$
The area of this upper bound is of order 
$(n^{1-\delta}/N) (\ml N/n)=\ml/n^\delta$ and tends to $0$ at an algebraic
rate in $\Lambda$ provided $\beta$
is less than $\delta$. Therefore, Lemma \Poisson\ implies
$$
  \max_{\Lambda\leq n\leq\Lambda^{1+\epsilon}} \sharp\calB_{1,n,N}=O_P(1)
$$
as $\Lambda$ tends to infinity.

Considering the bad points in $T_{5,n,N}^+$, we have, since $(g_n)$ is
equivalent to a nondecreasing sequence,
$$\displaylines{
  \Bigl|\fnn \sum_{i\in\calB_{1,n,N}}g_{n-NV_i} W_i^{-1/\alpha}\Bigr|
  \hfill\cr\hfill
  \eqalign{
  {}\leq{}&\fnn g_{n^{1-\delta}} \max_{i\in\RTwo} W_i^{-1/\alpha}
           \sharp\calB_{1,n,N}  \quad\equa{approxSixA}\cr
  {}={}   & \fnn g_n n^{-\delta(\gamma-1)}
            \Bigl({n\xi_N\over N}\Bigr)^{1/\alpha} O_P(1) \, . \cr}%
  \cr}
$$
This upper bound is of order $F^\leftarrow(1-1/n) 
g_n n^{-\delta(\gamma-1)}\xi^{1/\alpha}(N)$;  since $n$ is at 
least $\Lambda$ and $\xi$ is slowly varying, it is of order at most
$F^\leftarrow(1-1/n)g_n n^{-\delta(\gamma-1)/2}$. Lemma \trick\ then implies
$$\displaylines{
  \liml{\rm P}\Bigl\{\,\existsnl \fnn
  \hfill\cr\hfill\sum_{i\in\calB_{1,n,N}} g_{n-NV_i}W_i^{-1/\alpha}
  >T\g n/\kl\,\Bigr\} = 0 \, .
  \quad\equa{approxSixB}\cr}
$$
Referring to $T_{6,n,N}^+$,
$$\displaylines{\quad
  \fnn g_n\sum_{i\in \calB_{1,n,N}} \Bigl(1-{N\over n}V_i\Bigr)_+^{\gamma-1}
  W_i^{-1/\alpha}
  \hfill\cr\hfill
  {}\leq \fnn g_n n^{-\delta(\gamma-1)} \max_{i\in\RTwo}W_i^{-1/\alpha}
  \sharp\calB_{1,n,N} \, ,
  \quad\cr}
$$
which is the same bound as in \approxSixA. Therefore, the analogue of
\approxSixB\ holds when substituting $g_n (1-NV_i/n)_+^{\gamma-1}$ 
for $g_{n-NV_i}$.

Dealing with the bad points in $\calB_{2,n,N}$ is easy because if $V_i>n/N$
then $g_{n-NV_i}$ and $(1-NV_i/n)_+^{\gamma-1}$ vanish. So those points do 
not contribute to $T_{5,n,N}^+$ and $T_{6,n,N}^+$. 

On the part of $T_{5,n,N}^+$ made by the good points, assumption \hypMg\
and Lemma \cardinalityRegion\ yield, on the range of $n$ between $\Lambda$
and $\Lambda^{1+\epsilon}$,
$$\displaylines{
  \sum_{i\in\RTwo\setminus(\calB_{1,n,N}\cup\calB_{2,n,N})}
  \Bigl|g_{n-NV_i}-g_n\Bigl(1-{NV_i\over n}\Bigr)^{\gamma-1}\Bigr| 
  W_i^{-1/\alpha}
  \hfill\cr\hfill
  \eqalign{
  {}\leq{}&n^{-\delta} g_n\max_{i\in\RTwo} W_i^{-1/\alpha}
           \sharp\RTwo \cr
  {}\leq{}&\Lambda^{-\delta} g_n
           \Bigl({n\xi_N\over N}\Bigr)^{1/\alpha} \ml\log N O_P(1)\, .\cr
  }
  \qquad\cr}
$$
The bound obtained for the error in the 
approximation in $T_{5,n,N}^+$ for the good points is then
$$
  \fn g_n \xi_N^{1/\alpha} \Lambda^{-\delta}
  \ml \log N O_P(1) \, .
$$
If $\beta$ is less than $\delta$ then 
$\xi_N^{1/\alpha}\Lambda^{-\delta}\ml\log N$
tends to $0$ as $\Lambda$ tends to infinity and  Lemma \trick\ 
yields the conclusion.\hfill\qed

\bigskip

We now replace $\RTwo$ in $T_{6,n,N}^+$ by
$$
  \RThree=\Bigl\{\, i\,:\, W_i\leq \ml {N\over n} \, ,\, V_i\leq {n\over N}\,
          \Bigr\} \, ,
$$
defining
$$
  T_{7,n,N}^+
  = \fnn g_n \sum_{i\in\RThree}\kvi W_i^{-1/\alpha} \, .
$$

\Lemma{\label{approxSeven}
  If $\beta$ is less than $(\gamma-1)/2$,
  $$
    \liml\Prob{\existsnl |T_{6,n,N}^+-T_{7,n,N}^+|>T{\g n\over\kl}}
    =0 \, .
  $$
}

\Proof If $i$ belongs to $\RTwo\bigtriangleup\RThree$, then either
$$
  V_i\leq {n\over N} \leq G_N(V_i)
  \qquad\hbox{or}\qquad
  G_N(V_i)\leq {n\over N} \leq V_i \, .
$$
In the latter case, $(1-NV_i/n)_+$ vanishes and so those points do not 
contribute to $T_{6,n,N}^+$ and $T_{7,n,N}^+$. In the former case, arguing as
in the proof of Lemma \approxFive, 
$$\eqalign{
  V_i \leq {n\over N}\leq G_N(V_i)
  &{}\leq V_i+{V_i^{(1/2)-\eta}\over\sqrt N} O_P(1) \cr
  &{}\leq V_i +{n^{(1/2)-\eta}\over N} N^\eta O_P(1) \, , \cr
  }
$$
the last inequality coming from the first one and the $O_P(1)$-term being
uniform over $i$ in $\RTwo\bigtriangleup\RThree$ and $n$ between $\Lambda$
and $\Lambda^{1+\epsilon}$. In particular,
$$
  \Bigl| 1-{N\over n} V_i\Bigr|
  \leq n^{-(1/2)-\eta} N^\eta O_P(1) \, . 
$$
Consequently, using lemma \cardinalityRegion\ with $\RThree$ substituted
for $\RTwo$,
$$\displaylines{\qquad
  |T_{6,n,N}^+-T_{7,n,N}^+|
  \hfill\cr\hfill
  \eqalign{
    {}\leq{}& \fnn g_n (  n^{-(1/2)-\eta} N^\eta)^{\gamma-1} 
              \Bigl({n\over N}\xi_N\Bigr)^{1/\alpha} \cr
            &\hskip 120pt {}\times
              \sharp (\RTwo\bigtriangleup\RThree) O_P(1) \cr
    {}\leq{}& \fn g_n \xi_N^{2/\alpha} 
              (n^{-(1/2)-\eta} N^\eta)^{\gamma-1} \ml \log N O_P(1) \, .\cr}
  \cr}
$$
Note that $n^{-(1/2)-\eta}N^\eta$ is at 
most $\Lambda^{-(1/2)-\eta+(1+2\epsilon)\eta}$. 
Provided $\eta$ is small enough,
$$
  \Bigl( -{1\over 2} -\eta+ (1+2\epsilon)\eta\Bigr) (\gamma-1)+\beta
$$
is negative. We apply Lemma \trick\ to 
conclude.\hfill\qed

\bigskip

Having approximated all these $T_{i,n,N}^+$, we need to consider their
expected values, a much easier task. This requires us to have an estimate on
how close $\g n$ is to $ng_n/\gamma$ and how large $\mu_{n,\Lambda}^+$ is. We establish
those estimates in the next two lemmas.

\Lemma{\label{gngZero}
  If \hypMg\ holds, then, as $n$ tends to infinity,
  $$
    {\gamma\g n\over ng_n}
    = 1+o(n^{-\delta})\, .
  $$
}

\Proof We write $\gamma \g n/n g_n$ as 
$$
  {\gamma\over n}\sum_{0\leq i< n^{1-\delta}} {g_i\over g_n}
  + {\gamma\over n} \sum_{n^{1-\delta}\leq i< n} {g_i\over g_n} \, .
  \eqno{\equa{gngZeroNA}}
$$
Since $(g_n)$ is asymptotically equivalent to a monotone sequence, 
the first term in \gngZeroNA\ is at 
most $O(1)n^{-\delta}g_{n^{1-\delta}}/g_n=o(n^{-\delta})$.
Using \hypMg, the second term is at most
$$\displaylines{\qquad
    {\gamma\over n} \sum_{n^{1-\delta}< i\leq n} 
    \Bigl(\Bigl({i\over n}\Bigr)^{\gamma-1}+o(n^{-\delta})\Bigr)
  \hfill\cr\hfill
  \eqalign{
    {}={}&{1\over n^\gamma}\Bigl(\int_{n^{1-\delta}}^n \gamma x^{\gamma-1}\d x
          +O(1)\Bigr) + o(n^{-\delta}) \cr
    {}={}&1+O(n^{-\gamma\delta}) +o(n^{-\delta}) \, . \cr
  }
  \qquad\cr}
$$
The result follows since $\gamma$ exceeds $1$.\hfill\qed

\bigskip

The proof of Lemma \mun\ shows that if $\beta$ is small enough,
$$
  \limsupn \sup_{\Lambda\leq n\leq \Lambda^{1+\epsilon}}
  \Bigl| k(n)\mu_{n,\Lambda}^+-{\alpha\over\alpha-1}\ml^{1-1/\alpha}\Bigr|
  <\infty \, .
  \eqno{\equa{munA}}
$$

\bigskip

\Lemma{\label{approxExpectation}
  If $\beta<\delta\alpha/(\alpha-1)$, there exists a positive $T$ such that
  for any $\Lambda$ large enough and any $n$ 
  in $(\Lambda,\Lambda^{1+\epsilon})$,
  $$
    |\E T_{1,n,N}^+-\E T_{7,n,N}^+|
    \leq T {\g n\over k(\Lambda)} \, .
  $$
}

\Proof Lemmas \gngZero\ and \munA\ yield
$$\eqalignno{
  \E T_{1,n,N}^+- {ng_n\over\gamma}\mu_{n,\Lambda}^+
  &{}= n g_n\mu_{n,\Lambda}^+ o(n^{-\delta})  &\equa{approxExpectationA}\cr
  &{}= g_n \fn  \ml^{1-1/\alpha} n^{-\delta} o(1) \, . \cr
  }
$$
The sequence $\ml^{1-1/\alpha} \Lambda^{-\delta}$ is regularly
varying of index $\beta (1-1/\alpha)-\delta$, which is negative
provided $\beta$ is small enough. Thus, Lemma \trick\ shows that we can
replace $\E T_{1,n,N}^+$ by $n g_n \mu_{n,\Lambda}^+/\gamma$ to prove the current
lemma.

The calculation of $\E T_{7,n,N}^+$ can be done by using that 
$W_i$ has a gamma distribution with parameter $i$, but an easier argument
will show after our next lemma that
$$
  \E T_{7,n,N}^+
  = \fn {g_n\over\gamma}{\alpha\over\alpha-1} \ml^{1-1/\alpha} \, .
  \eqno{\equa{approxExpectationB}}
$$
Thus, it suffices to show that for any $T$ large enough,
$$
  \fn {g_n\over\gamma} \Bigl( k(n)\mu_{n,\Lambda}^+-{\alpha\over\alpha-1}
  \ml^{1-1/\alpha}\Bigr)
  \leq T \g n/\kl \, .
$$
This follows from \munA\ and Lemma \trick.\hfill\qed

\bigskip

Combining Lemmas \approxTwo\hskip 2pt--\hskip1.75pt\approxSeven\ and 
\approxExpectation, we see that
in order to prove Proposition \stepFour, it suffices to show that
$$
  \limsupL\prob{ \existsnl |T_{7,n,N}^+-\E T_{7,n,N}^+|
  \geq T\g n/\kl} = 0 \, .
  \eqno{\equa{stepFourA}}
$$

\bigskip

Our next lemma will allow us to represent $T_{7,n,N}^+$ as an integral
of a Poisson process over the quadrant $[\,0,\infty)^2$ and prove a valuable
scaling property. 

\Lemma{\label{poissonExtension}
  We can construct a Poisson process $\Pi$ on $[\,0,\infty)^2$ with mean
  intensity the Lebesgue measure, such that the point process obtained
  by restricting $\Pi$ to $[\,0,1\,]\times[\,0,\infty)$ coincides with
  $\sum_{i\geq 1} \delta_{(V_i,W_i)}$.
}

\bigskip

\Proof Let $N'$ be a homogenous and unit intensity Poisson random measure on 
$(1,\infty)\times [\,0,\infty)$, independent of the sequence $(V_i,W_i)$. 
Define $N$ as $N'+\sum_{i\geq 1}\delta_{(V_i,W_i)}$. \hfill\qed

\bigskip

We then rewrite $T_{7,n,N}^+$ as
$$\displaylines{\qquad
  \fn g_n \Bigl({N\over n}\Bigr)^{(1/\alpha)+\gamma-1}\sum_{i\geq 1}
  \Bigl({n\over N}-V_i\Bigr)_+^{\gamma-1} W_i^{-1/\alpha}
  \hfill\cr\hfill
  \One\Bigl\{ {n\over N}W_i\leq \ml\,\Bigr\} \, .\quad\cr}
$$
In this expression, considering only the sum over $i$ and
thinking of $n/N$ as a continuous variable $t$, we are led to introduce
the process
$$
  \Upsilon_\Lambda(t)
  =\sum_{i\geq 1} (t-V_i)_+^{\gamma-1} W_i^{-1/\alpha}\One\{\, tW_i\leq\ml\,\}
$$
indexed by $t$ in $[\,0,1\,]$.
Given Lemma \poissonExtension, we can extend $\Upsilon_\Lambda$ to a process
over the nonnegative half-line
$$
  \Upsilon_\Lambda(t)
  =\int(t-v)_+^{\gamma-1} w^{-1/\alpha}\One\{\, tw\leq\ml\,\}\d \Pi(v,w)
  \, .
$$
We then have
$$
  T_{7,n,N}^+=\fn g_n\Bigl({N\over n}\Bigr)^{(1/\alpha)+\gamma-1}
  \Upsilon_\Lambda\Bigl({n\over N}\Bigr) \, . 
  \eqno{\equa{stepFourAa}}
$$
In particular, $\E T_{7,n,N}^+$ is indeed given by \approxExpectationB\ since
the intensity of $\Pi$ being the Lebesgue measure,
$$\eqalign{
  \E\Upsilon_\Lambda(t)
  &{}= \int (t-v)_+^{\gamma-1} w^{-1/\alpha} \One\{\, tw\leq \ml\,\}\d v\d w 
    \cr
  &{}={t^{\gamma-1+1/\alpha}\over\gamma}\ml^{(\alpha-1)/\alpha}
    {\alpha\over\alpha-1} \, .\cr
  }
$$

Given \stepFourAa, in order to prove \stepFourA\ it suffices to show that
$$\displaylines{
  \limT\limsupL \Prob{\existsnl \fn g_n
  \Bigl({N\over n}\Bigr)^{(1/\alpha)+\gamma-1} 
  \hfill\cr\hfill
  |\Upsilon_\Lambda(n/N)-\E \Upsilon_\Lambda (n/N)| > T \g n/\kl}
  = 0 \, .
  \quad\equa{stepFourB}\cr}
$$

The scaling property we alluded to is the following.

\Lemma{\label{upsilonScaling}
  For any positive real number $\lambda$, the processes
  $\Upsilon_\Lambda(\lambda\, \cdot\,)$ and 
  $\lambda^{\gamma-1+(1/\alpha)}\Upsilon_\Lambda$ have the same distribution.
}

\bigskip

\Proof We rewrite $\Upsilon_\Lambda(\lambda t)$ as
\hfuzz=1pt
$$\displaylines{\quad
    \int(\lambda t-v)_+^{\gamma-1} w^{-1/\alpha} 
    \One\{\, \lambda t w\leq\ml\,\} \d \Pi(v,w) 
  \hfill\cr\hfill
    =\lambda^{\gamma-1+1/\alpha} \int 
    \Bigl( t-{v\over\lambda}\Bigr)_+^{\gamma-1}(\lambda w)^{-1/\alpha}
    \One\{\, t \lambda w\leq \ml\,\} \d \Pi(v,w) \, . 
  \cr}
$$
\hfuzz=0pt
The image of the Poisson random measure $\Pi$ by the map $(v,w)\mapsto 
(v/\lambda,\lambda w)$ is a Poisson random measure of intensity the Lebesgue
measure, proving the lemma.\hfill\qed

\bigskip

In what follows we will use the following terminology. We say that a 
sequence of either functions or random variables, $(f_n)$, converges to $f$ in
$\LTwo(\mu)$-norm if $\limn \int (f_n-f)^2\d\mu$ converges to $0$ as
$n$ tends to infinity. In our setting, $(f_n)$ and $f$ may not be in
$\LTwo(\mu)$ but $f_n-f$ is. If $\mu$ is the underlying probability ${\rm P}$,
we will write $\LTwo$ for $\LTwo({\rm P})$; in that case, convergence
of the sequence of random variables $(f_n)$ to $f$ in $\LTwo$-norm means 
that $\limn\E (f_n-f)^2=0$, again, even though $f_n$ and $f$ may not have 
finite variance but $f_n-f$ does.

Similarly, we will write that $f=g$ in $\LTwo$-norm, to mean $\E (f-g)^2=0$,
even though $f$ and $g$ may not be square integrable.

\bigskip

Writing $n/N$ as $(\Lambda/N)(n/\Lambda)$, Lemma \upsilonScaling\ shows that,
as a process indexed now by $n$ in $(\Lambda,\Lambda^{1+\epsilon})$,
$$
  \Upsilon_\Lambda \Bigl( {n\over N}\Bigr)
  \eqd \Bigl({\Lambda\over N}\Bigr)^{\gamma-1+(1/\alpha)} 
  \Upsilon_\Lambda\Bigl({n\over\Lambda}\Bigr) \, . 
$$
Therefore, to prove \stepFourB\ it suffices to show that
$$\displaylines{
  \limT\limsupL {\rm P}\Bigl\{\, \existsnl \fn g_n
  \Bigl({\Lambda\over n}\Bigr)^{(1/\alpha)+\gamma-1} 
  \hfill\cr\hfill
  \Bigl|\Upsilon_\Lambda
  \Bigl({n\over\Lambda}\Bigr)-\E\Upsilon\Bigl({n\over\Lambda}\Bigr)\Bigr|
  > T \g n/\kl\,\Bigr\} = 0 \, ,
  \cr}
$$
or, equivalently, that
$$\displaylines{\quad
    \limT\limsupL {\rm P}\Bigl\{\,\existsnl \Bigl|\Upsilon_\Lambda
    \Bigl({n\over\Lambda}\Bigr)-
    \E \Upsilon_\Lambda\Bigl({n\over\Lambda}\Bigr)\Bigr| 
  \hfill\cr\hfill
    {}>T 
    \Bigl({n\over\Lambda}\Bigr)^{(1/\alpha)+\gamma-1} {k(n)\over k(\Lambda)}
    \,\Bigr\} = 0 \, .
  \cr}
$$
Let $\eta$ be a positive real number less than $(\alpha-1)/2\alpha$
and $1-(\alpha/2)$. Using
Potter's bound, $k(n)/k(\Lambda)\gsim (n/\Lambda)^{1-(1/\alpha)-\eta}$
uniformly in $n$ between $\Lambda$ and $\Lambda^{1+\epsilon}$ and as $\Lambda$
tends to infinity. Thus to prove \stepFourB\ it suffices to show that
$$\displaylines{\quad
    \limT\limsupL\Prob{\existsnl 
  \hfill\cr\hfill
    \Bigl|\Upsilon_\Lambda\Bigl({n\over \Lambda}\Bigr)
    -\E \Upsilon_\Lambda\Bigl({n\over \Lambda}\Bigr)\Bigr| > T 
    \Bigl({n\over\Lambda}\Bigr)^{\gamma-\eta}} = 0 \, .
  \quad\equa{stepFourC}\cr}
$$

Setting 
$$
  f_{\Lambda,t}(v,w)
  =(t-v)_+^{\gamma-1} w^{-1/\alpha} \One\{\, w\leq m_\Lambda/t\,\}\, ,
$$
we rewrite the centered version of $\Upsilon_\Lambda$ as a compensated
Poisson integral,
$$
  (\Upsilon_\Lambda-\E \Upsilon_\Lambda)(t)
  = \int f_{\Lambda,t}(v,w) \d (\Pi-\E \Pi)(v,w) \, .
$$
As $\Lambda$ tends to infinity, the function $f_{\Lambda,t}$ converges 
pointwise to
$$
  f_t(v,w)=(t-v)_+^{\gamma-1} w^{-1/\alpha} \, .
$$
This convergence holds in $\LTwo(\d v\d w)$-norm since
$$\displaylines{\qquad
  \int (f_{\Lambda,t}-f_t)^2 (v,w)\d v\d w
  \hfill\cr\hfill
  \eqalign{
    {}={}&\int (t-v)_+^{2(\gamma-1)} w^{-2/\alpha} \One\{\, w\geq \ml /t\,\} 
          \d v \d w \cr
    {}={}&t^{2(\gamma-1+1/\alpha)} \ml^{-(2/\alpha)+1} 
          {\alpha\over(2-\alpha)(2\gamma-1)} \, . 
          \qquad\equa{stepFourD}\cr}\cr}
$$
Recall that the compensated Poisson integral induces an isometry in the sense
that for any function $f$ in $\LTwo(\d v \d w)$
$$
  \E \Bigl(\int f \d (\Pi-\E \Pi)\Bigr)^2 = \int f^2(v,w) \d v \d w \, .
$$
It then follows from \stepFourD\ that 
$\liml\Upsilon_\Lambda-\E\Upsilon_\Lambda$
exists pointwise in $\LTwo$-norm and is the compensated Poisson integral
$$
  \Upsilon_0(t)=\int f_t \d (\Pi-\E\Pi) \, .
$$
Our next lemma implies that we can replace 
$\Upsilon_\Lambda-\E\Upsilon_\Lambda$ by its limit $\Upsilon_0$ in \stepFour.

\Lemma{\label{approxUpsilon} Provided $\beta$ is positive, we have
  $$\displaylines{
    \limT\limsupL \Prob{\existsnl \Bigl| 
    \int (f_{\Lambda,n/\Lambda}-f_{n/\Lambda})\d (\Pi-\E\Pi)
    \Bigr| 
    \hfill\cr\hfill
    {}\geq T \Bigl({n\over \Lambda}\Bigr)^{\gamma-\eta} } = 0 \, .\cr}
  $$
}

\Proof Since 
$$
  (f_t-f_{\Lambda,t})(v,w)
  = (t-v)_+^{\gamma-1} w^{-1/\alpha}\One\{\, w\geq \ml/t\,\}
$$
is in ${\rm L}^p(\d v \d w)$ for any $p$ 
greater than $\alpha$, let
$$
  M_p=\E\Bigl( \int (f_{\Lambda,t}-f_t)\d (\Pi-\E\Pi) \Bigr)^p \, .
$$
In what follows we restrict $p$ to be an integer. Using Privault's (2009,
equation 9; or 2010, equation 2.9) moment identity (see also Bassan and Bona, 
1990), we have for $p$ at least
$2$,
$$
  M_p=\sum_{0\leq k\leq p-2} {p-1\choose k} \int (f_{\Lambda,t}-f_t)^{p-k}(v,w)
  \d v \d w M_k \, .
  \eqno{\equa{approxUpsilonA}}
$$
We now prove by induction that for some constant $c_p$,
$$
  |M_p|\leq c_p t^{p(\gamma-1+1/\alpha)} \ml^{-(p/\alpha)+\lfloor p/2\rfloor}
  \, .
  \eqno{\equa{approxUpsilonB}}
$$
Indeed, $M_1$ vanishes and, as shown in \stepFourD,
$$
  M_2 =c_2 t^{2(\gamma-1+1/\alpha)}\ml^{-(2/\alpha)+1} \, .
$$
Assume that for any $k$ less than $p$,
$$
  |M_k|\leq c_k t^{k(\gamma-1+1/\alpha)} \ml^{-(k/\alpha)+\lfloor k/2\rfloor}
  \, . 
$$
Since
$$
  \int (f_t-f_{\Lambda,t})^{p-k}(v,w)\d v \d w
  = c t^{(p-k)(\gamma-1+1/\alpha)} \ml^{1-((p-k)/\alpha)} \,  ,
$$
equality \approxUpsilonA\ and the induction hypothesis imply
$$\displaylines{\qquad
  |M_p| 
  \leq c \sum_{0\leq k\leq p-2} t^{(p-k)(\gamma-1+1/\alpha)}
       \ml^{1-(p-k)/\alpha} 
  \hfill\cr\hfill
    {}\times t^{k(\gamma-1+1/\alpha)} 
       \ml^{-(k/\alpha)+\lfloor k/2\rfloor} \, . 
  \qquad\cr}
$$
In this sum, bounding $\ml^{-(k/\alpha)+\lfloor k/2\rfloor}$ 
by $\ml^{-(k/\alpha)+\lfloor (p-2)/2\rfloor}$ yields
$$
  |M_p| \leq c_p t^{p(\gamma-1+1/\alpha)} 
  \ml^{1-(p/\alpha)+\lfloor (p-2)/2\rfloor} \, , 
$$
which is \approxUpsilonB.

We then take $p$ to be an even integer. Applying Markov's inequality for 
any $n$ between $\Lambda$ and $\Lambda^{1+\epsilon}$ and using \approxUpsilonB,
$$\displaylines{\qquad
  \Prob{ \Bigl| \int (f_{\Lambda,n/\Lambda}-f_{n/\Lambda})\d (\Pi-\E\Pi)\Bigr|
  \geq T \Bigl({n\over\Lambda}\Bigr)^{\gamma-\eta} }
  \hfill\cr\hfill
  \eqalign{
    {}\leq{}& T^{-p} \Bigl({\Lambda\over n}\Bigr)^{(\gamma-\eta) p} c_p
              \Bigl({n\over\Lambda}\Bigr) ^{p(\gamma-1+1/\alpha)}
              \ml^{-(p/\alpha)+\lfloor p/2\rfloor} \cr
    {}\leq{}& c_p T^{-p} \Bigl({\Lambda\over n}\Bigr)^{p(1-(1/\alpha)-\eta)}
              \ml^{-p((1/\alpha)-(1/2))} \, . \cr
  }\qquad\cr}
$$
Applying Bonferroni's inequality, the probability involved in the lemma is at
most
$$
  c_p T^{-p} \Lambda^{p(1-1/\alpha-\eta)} \ml^{-p((1/\alpha)-(1/2))}
  \sum_{\Lambda\leq n\leq \Lambda^{1+\epsilon}} 
  n^{-p(1-(1/\alpha)-\eta)} \, .
$$
Taking $p$ larger than $1/\bigl(1-(1/\alpha)-\eta\bigr)$, this bound is of 
order
$$\displaylines{\qquad
    c T^{-p} \Lambda^{p(1-(1/\alpha))} \ml^{-p(1-(1/\alpha)-(1/2))}
    \Lambda^{1-p(1-(1/\alpha))}
  \hfill\cr\hfill
    {}=c T^{-p}\Lambda\ml^{-p((1/\alpha)-(1/2))} \, .\qquad\cr}
$$
This bound is regularly varying in $\Lambda$ of index
$$
  1-\beta p\Bigl( {1\over\alpha}-{1\over 2}\Bigr) \, .
$$
Thus, taking $p$ larger than $2\alpha/\bigl(\beta (2-\alpha)\bigr)$, it tends
to $0$ as $\Lambda$ tends to infinity, proving the lemma.\hfill\qed

\bigskip

Given Lemma \approxUpsilon, we see that to prove \stepFourC\ it suffices to
show that 
$$
  \limT\limsupL \Prob{ \existsnl |\Upsilon_0(n/\Lambda)|
  >T (n/\Lambda)^{\gamma-\eta}} = 0 \, . 
  \eqno{\equa{stepFourE}}
$$

Our next step is to identify the process $\Upsilon_0$ as a fractional integral
of a spectrally positive L\'evy stable process.

Consider the spectrally positive centered L\'evy stable process given by
its It\^o representation
$$
  L_0^+(t)=\int \One_{[0,t]}(v) w^{-1/\alpha} \d (\Pi-\E\Pi)(v,w) \, .
$$ 
This process is the pointwise limit in $\LTwo$-norm of
$$
  L_{0,\epsilon}^+(t)
  =\int \One_{[0,t]}(v) w^{-1/\alpha} \One_{[0,1/\epsilon]}(w)
  \d (\Pi-\E\Pi)(v,w)
$$
as $\epsilon$ tends to $0$. Moreover, defining
$$\eqalign{
  L_{0,\epsilon}^{+(\gamma-1)}(t)
  &{}=\int \gamma(t-v)_+^{\gamma-1}\d L_{0,\epsilon}^+(v) \cr
  &{}=\sum_{i\geq 1} \Bigl( 
    \gamma(t-V_i)_+^{\gamma-1}W_i^{-1/\alpha}\One_{[0,1/\epsilon]}(W_i) \cr
  &\hskip 120pt
    {}- t^\gamma \E W_i^{-1/\alpha} \One_{[0,1/\epsilon]}(W_i)\Bigr) \cr
  &{}=\int\gamma(t-v)_+^{\gamma-1} w^{-1/\alpha} \One_{[0,1/\epsilon]}(w)
    \d (\Pi-\E\Pi)(v,w) \, , \cr
  }
$$
we see that $L_{0,\epsilon}^{+(\gamma-1)}$ converges pointwise in $\LTwo$-norm
to $\gamma\Upsilon_0$ as $\epsilon$ tends to $0$. It follows that 
pointwise in $\LTwo$-norm,
$$
  \Upsilon_0^+(t)=\int (t-v)_+^{\gamma-1} \d L_0^+(v) \, .
  \eqno{\equa{stepFourF}}
$$
Lemma 3.1.7 in \gFLevy\ implies that the right hand side of \stepFourF\ is 
almost surely continuous. Considering the fractional integral
$$
  L_0^{+(\gamma-1)}(t)=\int \gamma (t-v)_+^{\gamma-1} \d L_0^+(v) \, ,
$$
and since $\Upsilon_0^+(n/\Lambda)$ and $L_0^{+(\gamma-1)}(n/\Lambda)/\gamma$ 
coincide in $\LTwo$-norm for any integer $n$ between $\Lambda$ 
and $\Lambda^{1+\epsilon}$,
we see that in order to prove $\stepFourE$ and therefore Proposition \stepFour,
it suffices to show the following.

\Lemma{\label{limitBounded}
  For any positive $\eta$ sufficiently small,
  $$
    \limT\prob{\exists t\geq 1\,:\, |L_0^{+(\gamma-1)}(t)|\geq 
    T t^{\gamma-\eta}}= 0 \, .
  $$
}

\Proof Note that $L_0^+$ vanishes at $0$. For any nonnegative $t$, the 
function 
$(t-\Id)_+^{\gamma-1}$ is deterministic, differentiable 
on $[\,0,t\,]$, so that its quadratic covariation with $L_0^+$ vanishes
on $[\,0,t\,]$. We then integrate by parts the integral 
defining $L_0^{+(\gamma-1)}$ (see Protter, 1992, chapter 2.6, Corollary 2)
as
$$
  L_0^{+(\gamma-1)}(t)
  =\gamma (\gamma-1)\int L^+_0(v) (t-v)_+^{\gamma-2} \d v.
$$
It follows from Pruitt (1981) that 
there exists a constant $c$ and a random $v_0$ such that
$|L^+_0(v)|\leq c v^{(1/\alpha)+\eta}$ for any $v$ at least $v_0$. 
This implies that for $t$ at least $v_0$,
$$ 
  \int_{v_0}^t |L^+_0(v)| (t-v)_+^{\gamma-2} \d v
  \leq c t^{\gamma-1+(1/\alpha)+\eta} \, .
$$
Since 
$$
  \int_0^{v_0} |L^+_0(v)|(t-v)_+^{\gamma-2}\d v
  \leq c \sup_{0\leq v\leq v_0} |L_0^+(v)| t^{\gamma-1}
$$
and $-1+(1/\alpha)+2\eta$ is negative for $\eta$ small enough, we have
$$
  \limt t^{-\gamma+\eta}L_0^{+(\gamma-1)}(t)=0
$$
almost surely. The lemma follows.\hfill\qed

\bigskip


\subsection{Concluding the proof}
Having completed the proof of Proposition \stepFour, we can complete that of
Theorem \GammaGTOne.

First, we settle the assertion that
$$
  \sup_{t\geq 0} L_0^{(\gamma-1)}(t)-t^\gamma 
  \hbox{ is almost surely finite.}
  \eqno{\equa{tight}}
$$
Let $\widetilde L_0^{+(\gamma-1)}$ be an independent copy 
of $L_0^{+(\gamma-1)}$. It
is shown in section 3.3 of \gFLevy\ that $L_0^{(\gamma-1)}$ has the same
distribution as $p^{1/\alpha}L_0^{+(\gamma-1)}-q^{1/\alpha}
\widetilde L_0^{+(\gamma-1)}$.
Consequently, \tight\ follows from Lemma \limitBounded.

To complete the proof of Theorem \GammaGTOne, we set
$$
  T_n^-=\sum_{0\leq i<n} g_i X_{n-i}\One\{\, X_{n-i}\leq a_{n,\Lambda}\,\} \, .
$$
Note that
$$
  T_n^-= -\sum_{0\leq i<n} g_i (-X_{n-i})\One\{\, X_{n-i}\geq -a_{n,\Lambda}\,\} \, .
$$
If assumption \hypRVRateLT\ hold,
we substitute $\M_{-1}F$ for $F$ in steps two, three and four to obtain 
that
$$
  \limT\limsupL\prob{\exists n\geq\Lambda\,:\,
  |T_n^--\E T_n^-|>T\g n/\kl} = 0 \, .
$$
This, combined with Propositions \middleNeglect, \stepTwo, \middleNeglectThree\
and \stepFour\ give our heavy traffic approximation under \hypRVRateLT.

\bigskip

In order to prove Theorem \GammaGTOne\ under \tailDominance\ we will use a
coupling argument based on a decomposition of the innovations into a part
which is bounded from above and a part which is bounded from below. The
underlying idea is that the part bounded from above should not contribute
too much to the process $S_n$ reaching the boundary $a\g n$. 
To make this argument viable, note that combining Propositions \middleNeglect,
\stepTwo, \middleNeglectThree\ and \stepFour, Theorem \GammaGTOne\ holds
when $F$ is supported on some interval bounded away from $-\infty$.

The following lemma allows us to define the proper representation of the 
innovations. If $G$ is a distribution function, $G^\leftarrow(0+)$ and
$G^\leftarrow(1)$ are respectively the lower and upper end point of the 
support of the underlying probability measure.

\Lemma{\label{representation}
  Let $F$ be a distribution function on the real line, centered and such that
  $\oF$ is regularly varying. There exists two distribution functions $F_L$ and
  $F_U$ and a $r$ in $(0,1)$ such that
  
  \medskip

  \noindent (i) $F_U^\leftarrow(0+)>-\infty$ and $F_L^\leftarrow(1)<\infty$,

  \medskip

  \noindent (ii) $F_U$ and $F_L$ are centered,

  \medskip

  \noindent (iii) $F=rF_U+(1-r)F_L$.
}

\bigskip

\noindent
Given assertions (i) and (iii), we must have $\oF=r\oF_U$ 
and $\overline{\M_{-1}F}=(1-r)\overline{\M_{-1}F_L}$ ultimately. 
Thus $rF_U$ and $(1-r)F_L$ capture respectively the upper and lower 
tails of $F$.

\bigskip

\Proof Since $F$ has a mean, define $A(t)=\int_{(t,\infty)}x\d F(x)$ and 
set $\theta=A(0)/2$.
The function $A$ is nonincreasing and c\`adl\`ag on the nonnegative half-line.
Let 
$$
  t_1=\inf\{\, t\,:\, A(t)<\theta\,\} \, .
$$
If $A(t_1)<\theta$, then $F$ has a jump at $t_1$ and there
exists a positive $\tau_1$ at most $F(t_1)-F(t_1-)$ such 
that $A(t_1)+\tau_1 t_1=\theta$; otherwise we set $\tau_1=0$.

Similarly, define $B(t)=\int_{(t,0]}x\d F(x)$. The function $B$ is 
nondecreasing on the negative half-line. Since $F$ is centered, 
$\lim_{t\to-\infty} B(t)=-A(0)$. Let 
$$
  t_2=\sup\{\, t\,:\, B(t)<-\theta\,\} \, .
$$ 
If $B(t_2)>-\theta$ then $F$ has a jump at $t_2$ and we define $\tau_2$
such that $B(t_2-)+\tau_2 t_2=-\theta$; otherwise, we set $\tau_2=0$.

We define the measure $\mu_U$ by
$$
  {\d\mu_U\over \d F}=\One_{(t_2,0]}+\One_{(t_1,\infty)}+\tau_1\One_{\{t_1\}}
  + \tau_2\One_{\{ t_2\}} \, ,
$$
and set $F_U=\mu_U/\mu_U(\RR)$. By construction $F_U$ is centered. 
Set $r=\mu_U(\RR)$ and define $F_L=(F-rF_U)/(1-r)$.\hfill\qed

\bigskip

Our next lemma, valid since $\alpha$ is positive and less than $2$,
relates the moment generating function of $F_L$ to the
bound \tailDominance. It is convenient to define
$$
  D(t)={c\over 1-r} \oF (t\log t)\log t \, ,
$$
$r$ being as in Lemma \representation\ and
the constant $c$ being, for once, the same as in \tailDominance

\Lemma{\label{mgf}
  The following holds as $\lambda$ tends to $0$ from above,
  $$
    \int e^{\lambda x} \d F_L(x)
    \leq 1+D\Bigl({1\over\lambda}\Bigr) \int_0^\infty {1-e^{-u}\over u^\alpha}
    \d u \bigl( 1+o(1)\bigr) \, .
  $$
}

\Proof Since $F_L$ is centered, two integrations by parts yield
$$
  0=\int x \d F_L(x)
  = -\int_{-\infty}^0 F_L(x)\d x +\int_0^\infty \oF_L(x)\d x \, .
$$
This identity, an integration by parts and considering that $\oF_L$ vanishes
ultimately yield that the moment generating function of $F_L$ is
$$\displaylines{\quad
  \lambda\int e^{\lambda x} \oF_L(x) \d x 
  \hfill\cr\hfill
  \eqalign{
    {}={}& 1-\lambda\int_{-\infty}^0 e^{\lambda x} F_L(x) \d x 
           +\lambda\int_0^\infty e^{\lambda x} \oF_L(x) \d x \cr
    {}={}& 1+\lambda\int_{-\infty}^0 (1-e^{\lambda x}) F_L(x)\d x
           +\lambda \int_0^\infty (e^{\lambda x}-1)\oF_L(x) \d x \, . \cr
  }\quad\cr
  }
$$
Let $t_0$ be a negative real number such that $(1-r)F_L$ coincides with $F$
on  $(-\infty,t_0\,]$ and and $\overline{M_{-1}F}\leq (1-r)D$ 
on $[\,t_0,\infty)$. Since $\oF_L$ vanishes ultimately,
$$
  \lambda \int_{t_0}^0 (1-e^{\lambda x})F_L(x)\d x
  +\lambda \int_0^\infty (e^{\lambda x}-1)\oF_L(x) \d x
  = O(\lambda^2)
$$
as $\lambda$ tends to $0$. Hence,
$$\eqalignno{
  \int e^{\lambda x}\d F_L(x)
  &{}\leq 1+\lambda \int_{-\infty}^{t_0} (1-e^{\lambda x})D(-x)\d x 
    + O(\lambda^2) \cr
  &{}= 1+\int_{-\infty}^{\lambda t_0}(1-e^y)D(-y/\lambda)\d y 
    +O(\lambda^2)\, . \qquad
  &\equa{mgfA}\cr
  }
$$
Since $D$ is regularly varying of index $-\alpha$ less than $-1$, the 
integral in \mgfA\ is asymptotically equivalent to
$$
  D\Bigl({1\over\lambda}\Bigr) \int_{-\infty}^0 {1-e^y\over (-y)^\alpha}\d y
  \, .
$$
Since $\alpha$ is between $1$ and $2$, $D(1/\lambda)\gg\lambda^2$ as $\lambda$
tends to $0$ and the lemma follows.\hfill\qed

\bigskip

To conclude the proof of Theorem \GammaGTOne, write $F=rF_U+(1-r)F_L$ for
the decomposition given in Lemma \representation. Let $(X_{U,i})$ and
$(X_{L,i})$ be two independent sequences of independent random variables
distributed respectively according to $F_U$ and $F_L$. Let $(B_i)$
be an independent sequence of independent random variables having a Bernoulli
distribution with parameter $p$. Set
$$
  X_i = B_i X_{U,i}+(1-B_i)X_{L,i} \, .
$$
By construction $(X_i)$ is a sequence of independent random variables
all having distribution $F$. Let $(S_n)$, $(S_{U,n})$ and $(S_{L,n})$
by the \hbox{$(g,F)$-{},} $(g,F_U)$-{} and $(g,F_L)$-processes based 
on $(X_i)$, $(X_{U,i})$ and $(X_{L,i})$ respectively. 
We have $S_n=S_{L,n}+S_{U,n}$. Hence, refering to \ThATightC,
$$\displaylines{\qquad
  \prob{\exists n\geq \Lambda \,:\, S_n>2 T\g n /\kl}
  \hfill\cr\noalign{\vskip 3pt}\hfill
  \leq 
  \prob{\exists n\geq \Lambda \,:\, S_{L,n}>T\g n /\kl}
  \qquad\qquad\cr\noalign{\vskip 3pt}\hfill
  {}+\prob{\exists n\geq \Lambda \,:\, S_{U,n}>T\g n /\kl}
  \quad\cr
}
$$
Using the part of Theorem \GammaGTOne\ that we proved already, more precisely
when the innovations are bounded from below,
$$
  \limT\limsupL \prob{\exists n\geq \Lambda \,:\, S_{U,n}>T\g n /\kl}
  = 0 \, .
$$
Hence, using a Bonferroni inequality, it suffices to show that
$$
  \limT\limsupL \sum_{n\geq\Lambda} \prob{S_{L,n}>T\g n /\kl}
  = 0 \, .
  \eqno{\equa{endA}}
$$
For this, let $\varphi$ be the moment generating function of $(1-B_i)X_{L,i}$
and let $\varphi_L$ be that of $X_{L,i}$. We have
$$
  \varphi(t)
  = \E e^{t(1-B_i)X_{L,i}}
  = r+(1-r)\varphi_L(t) \, .
  \eqno{\equa{endB}}
$$
The exponential form of Markov's inequality yields for any positive $\lambda$
$$
  \prob{ S_{L,n}>T\g n/\kl}
  \leq\exp\Bigl( -\lambda T {\g n\over\kl} +\sum_{0\leq i<n}
                 \log\varphi(\lambda g_i)\Bigr)\, .
  \eqno{\equa{endC}}
$$
In this inequality we take $\lambda=k(n)(\log n)/\g n$ and proceed to
estimate \endC.

Since $k$ is regularly varying of positive index, $k(\Lambda)\leq ck(n)$
if $n\geq\Lambda$ and $\Lambda$ is large enough. Consequently,
$\lambda T\g n/k(\Lambda)>c T$. 

Since \Karamata\ holds, our choice of $\lambda$ implies
$$
  \lambda\sim c{\log n\over F^\leftarrow(1-1/n) g_n} \, ;
$$
thus, since $(g_i)$ is asymptotically equivalent to a nondecreasing sequence, 
the bound
$$
  \max_{0\leq i<n}\lambda g_i \leq c {\log n\over F^\leftarrow(1-1/n)}
$$
shows that this maximum tends to $0$ as $\Lambda$ tends to infinity and $n$ is
at least $\Lambda$. Then, considering \endB\ and Lemma \mgf\ to 
bound $\varphi_L$, we obtain that for $\Lambda$ large enough, $n$ 
at least $\Lambda$ and any $i$ between $0$ and $n$,
$$\eqalign{
  \log\varphi(\lambda g_i)
  &{}\leq \log\Bigl( r+(1-r)\bigl(1+c D(1/\lambda g_i)\bigr)\Bigr)\cr
  &{}=\log\bigl(1+cD(1/\lambda g_i)\bigr) \cr
  &{}\leq c D(1/\lambda g_i) \, . \cr
  }
$$
Since $D$ is regularly varying of index $-\alpha$, this implies that $\endC$
is at most
$$
  \exp\Bigl( -Tc\log n +c n D\Bigl({F^\leftarrow(1-1/n)\over\log n}\Bigr)\Bigr)
  \, . 
  \eqno{\equa{endD}}
$$
We claim that
$$
  D\Bigl({F^\leftarrow(1-1/n)\over\log n}\Bigr)
  \leq c{\log n\over n}\, .
  \eqno{\equa{endE}}
$$
Indeed, setting $s=F^\leftarrow(1-1/n)/\log n$, we 
have $s\log s\sim c F^\leftarrow(1-1/n)$, which implies 
$\oF (s\log s)\sim c/n$. Moreover $\log s\sim c\log n$ as $n$ tends to 
infinity. This allows us to write \endE\ as $D(s)\leq c\log s \oF(s\log s)$
which is true by definition of $D$.

Combining \endD\ and \endE, we see that, if $T$ is large enough, \endC\ is
at most $n^{-2}$. This proves \endA\ and completes our proof 
of Theorem \GammaGTOne.

\bigskip


\subsection{On assumption \hypRVRate} To explain the role of assumption
\hypRVRate, we first explain that of $(m_n)$. The proof of Propositions
\middleNeglect, \middleNeglectThree\ and Lemma \approxUpsilon\ are similar: we 
evaluate some moment of high order and use Bonferroni's inequality. In 
order for the Bonferroni bound to tend to $0$, we need $(m_n)$ grow to
infinity at least at an algebraic rate.

In \fixedref{step 2} we do not rely on \hypRVRate\ but on \hypRVRateAlt. If
$(m_n)$ were allowed to be slowly varying, \fixedref{step 2} would sill hold.
The same is true for the approximations of the $T_{i,n,N}^+$ in 
\fixedref{step 4}, with the caveat that for Lemma \cardinalityRegion\ to hold,
we need to have $\liminf_{n\to\infty} m_n/\log n>0$, but this is essentially
implied by \hypMnUpsilon.

That $(m_n)$ has to grow at an algebraic rate makes \hypRVRateAlt\ equivalent
to \hypRVRate. So we see that the only reason \hypRVRate\ is needed is because
our crude estimate in the proofs of Propositions \middleNeglect, 
\middleNeglectThree\ and Lemma \approxUpsilon.

If we do not take those propositions and this lemma into account, the 
condition \hypMnUpsilon\ would allow $m_n=(\log n)^{1+\epsilon}$. To obtain
an even slower rate requires not to rely on Kiefer's theorem. It is
conceivable that a good description of the extremes 
$(U_i\One\{\, U_i>1-m_n/n\,\})_{1\leq i\leq n }$ as $n$ changes could yield
a better result with no other condition than \hypRVRate. This seems related to
the asymptotic behavior of $m_n$-records, but our attempt to devise a proof
in this direction failed.

Note that assumption \hypMg\ was used only in the proof of Lemma \approxSix\ 
and \approxExpectation. A close examination of the proofs shows that we could
replace \hypMg\ by the existence of a sequence $(\rho_n)$ such that
$$
  \limsup_{n\to\infty} \sup_{1/\rho_n\leq i/n\leq 1}\Bigl|{g_i\over g_n}
  -\Bigl({i\over n}\Bigr)^{\gamma-1}\Bigr| <\infty
$$
and $\rho_n\geq \xi_n^{1/\alpha}(m_n\vee\log n)$. So, any improvement
on the rate of growth of $(m_n)$ would translate into an improvement of
\hypMg\ as well.

\bigskip


\def\prevs{\the\sectionnumber .\the\snumber }
\def\preveq{(\the\sectionnumber .\the\equanumber)}

\section{Proof of Proposition \hypRVMeaning}
{\it (i)$\Rightarrow$(ii).} Write $F^\leftarrow(1-1/t)=t^{1/\alpha}
\ell(t)$ where $\ell$ is regularly varying. Note that $\ell$ is ultimately
positive. If (i) holds then, using $1/\lambda$ instead of $\lambda$,
$$
  \limt t^\kappa \sup_{t^{-\kappa}\leq\lambda \leq t^\kappa}
  \lambda^{1/\alpha} \Bigl|{\ell(\lambda t)\over\ell(t)}-1\Bigr| = 0 \, .
$$
In particular, for any fixed $\lambda$ greater than $1$, we have, for any
$t$ large enough,
$$
  \ell(t)(1-t^{-\kappa})
  \leq \ell(\lambda t)
  \leq \ell(t)(1+t^{-\kappa}) \, .
$$
Hence, for any positive integer $n$,
$$
  \ell(\lambda^{n-1}t)\bigl( 1-(\lambda^{n-1}t)^{-\kappa}\bigr)
  \leq \ell(\lambda^n t)
  \leq \ell(\lambda^{n-1}t)\bigl( 1+(\lambda^{n-1}t)^{-\kappa}\bigr) \, .
$$
By induction, this implies
$$
  \ell(t)\prod_{0\leq i<n}\bigl( 1-(\lambda^i t)^{-\kappa}\bigr)
  \leq \ell(\lambda^n t)
  \leq \ell(t)\prod_{0\leq i<n} \bigl(1+(\lambda^i t)^{-\kappa}\bigr) \, .
$$
Since $\ell$ is slowly varying,
$$
  \limn \sup_{s\in [\lambda^n t,\lambda^{n+1}t]}{\ell(s)\over\ell(\lambda^n t)}
  = 1\, .
$$
Therefore,
$$\eqalignno{
  \ell(t)\prod_{i\geq 0} \bigl( 1-(\lambda^it)^{-\kappa}\bigr)
  &{}\leq \liminf_{s\to\infty} \ell (s)\cr
  &{}\leq \limsup_{s\to\infty}\ell (s)
  \leq \ell(t)\prod_{i\geq 0} \bigl( 1+(\lambda^it)^{-\kappa}\bigr) \, .\qquad
  &\equa{rvA}\cr}
$$
Since 
$$
  \limt\prod_{i\geq 0}\bigl( 1-(\lambda^it)^{-\kappa}\bigr)
  =  \limt\prod_{i\geq 0}\bigl( 1+(\lambda^it)^{-\kappa}\bigr)
  =1 \, ,
$$
we obtain
$$
  \limsupt \ell(t)
  \leq\liminf_{s\to\infty}\ell(s)
  \leq\limsup_{s\to\infty}\ell(s)
  \leq\liminft\ell(t) \, ,
$$
proving that $\lim_{t\to\infty}\ell (t)$ exists. This limit is then 
positive due to \rvA.
Therefore, there exists a function $\delta(\cdot)$ which
tends to $0$ at infinity such that
$$
  F^\leftarrow(1-1/t)=ct^{1/\alpha}\bigl(1+\delta(t)\bigr) \, .
$$

We can then rewrite (i) as
$$
  \limt t^\kappa\sup_{t^{-\kappa}\leq\lambda\leq t^\kappa} \lambda^{1/\alpha}
  \Bigl| {1+\delta(\lambda t)\over 1+\delta(t)} -1\Bigr| = 0 \, .
  \eqno{\equa{rvB}}
$$
Since $\delta(\cdot)$ tends to $0$ at infinity, \rvB\ implies
$$
  \limt t^\kappa\sup_{t^{-\kappa}\leq\lambda\leq t^\kappa} \lambda^{1/\alpha}
  |\delta(\lambda t)-\delta(t)| = 0 \, .
$$
In particular, if $\epsilon$ is a fixed positive real number, for any $t$
large enough, considering $\lambda=t^\kappa$,
$$
  |\delta(t^{1+\kappa})-\delta(t)|\leq\epsilon t^{-\kappa(1+1/\alpha)} \, .
$$
Substituting $t^{(1+\kappa)^n}$ for $t$, we obtain
$$
  |\delta(t^{(1+\kappa)^{(n+1)}})-\delta(t^{(1+\kappa)^n})|
  \leq\epsilon t^{-\kappa(1+\kappa)^n(1+1/\alpha)} \, .
$$
Since $\delta(\cdot)$ tends to $0$ at infinity, summing all these inequalities
over $n$ nonnegative, we obtain
$$
  |\delta(t)|\leq \epsilon t^{-\kappa(1+1/\alpha)}
  +c t^{-\kappa(1+\kappa)(1+1/\alpha)} \, .
$$
Consequently, $\delta(t)=o(t^{-\kappa(1+1/\alpha)})$ and (ii) holds.

\smallskip

\noindent {\it (ii)$\Rightarrow$(i).} Whenever $t\wedge
(\lambda t)$ tends to infinity,
$$\eqalign{
  {F^\leftarrow(1-1/\lambda t)\over F^\leftarrow(1-1/t)}-\lambda^{1/\alpha}
  &{}=\lambda^{1/\alpha} 
    \Bigl({1+o\bigl((\lambda t)^{-\kappa(1+1/\alpha)}\bigr)\over
    1+o(t^{-\kappa(1+1/\alpha)})}-1\Bigr) \cr
  &{}=\lambda^{1/\alpha}\Bigl(o\bigl((\lambda t)^{-\kappa(1+1/\alpha)}\bigr)
    \vee o(t^{-\kappa(1+1/\alpha)})\Bigr) \, . \cr
  }
$$
In particular, since $\kappa$ is less than both $1$ and $2/(\alpha+1)$,
$$
  \sup_{t^{-\kappa}\leq\lambda\leq t^\kappa}
  \Bigl| {F^\leftarrow(1-1/\lambda t)\over F^\leftarrow(1-1/t)}
  -\lambda^{1/\alpha}\Bigr| = o(t^{-\kappa})
$$
and (i) holds.

\smallskip

\noindent{\it (ii)$\Rightarrow$(iii).} Consider the relation 
$x=F^\leftarrow(1-1/t)$. Given (ii), this means, as either $t$ or $x$ tend
to infinity,
$$
  x=ct^{1/\alpha}\bigl(1+o(t^{-\kappa(1+1/\alpha)})\bigr) \, .
  \eqno{\equa{rvD}}
$$
The proof is then an easy exercise in asymptotic analysis.
In particular, $x\sim c t^{1/\alpha}$ and $t\sim (x/c)^\alpha$. Write 
$t=(x/c)^\alpha(1+y)$ where $y$ tends to $0$ as $t$ tends to infinity. Using
\rvD\ we obtain
$$\eqalign{
  x&{}=x(1+y)^{1/\alpha}\bigl(1+o(x^{-\kappa(\alpha+1)})\bigr) \cr
   &{}=x+{1\over\alpha}xy + o(x^{1-\kappa(\alpha+1)}) + O(xy^2) \, . \cr
  }
$$
Therefore, $y= o(x^{-\kappa(\alpha+1)})+O(y^2)$. This 
forces $y=o(x^{-\kappa(\alpha+1)})$. Thus,
$$
  t=(x/c)^\alpha \bigl(1+o(x^{-\kappa(\alpha+1)})\bigr) \, .
  \eqno{\equa{rvE}}
$$
Since $x=F^\leftarrow(1-1/t)$ and $F$ is continuous, we have $\oF(x)=1/t$
as $t$ tends to infinity. Given \rvE, this implies
$$
  \oF(x)=(c/x)^\alpha\bigl(1+o(x^{-\kappa(\alpha+1)})\bigr) \, ,
$$
that is (iii).

\noindent{\it (iii)$\Rightarrow$(ii).} The proof is similar to that of its
converse.

\bigskip


\noindent {\bf References.}

\medskip

{\leftskip=\parindent \parindent=-\parindent
 \par

J.\ Akonom, Chr.\ Gouri\'eroux (1987). A functional central limit theorem
for fractional processes, discussion paper 8801, CEPREMAP, Paris.


Ph.\ Barbe, M.\ Broniatowski (1998). Note on functional large deviation
principle for fractional ARIMA processes, {\it Statistical Inference
for Stochastic Processes}, 1, 17--27.

Ph.\ Barbe, W.P.\ McCormick (2008a). Veraverbeke's theorem at large --
On the maximum of some processes with negative drift and heavy tail
innovations, {\tt arxiv:0802.3638}.


Ph.\ Barbe, W.P.\ McCormick (2008b). An extension of a logarithmic form of
Cram\'er's ruin theorem to some FARIMA and related processes,
{\tt arxiv:0811.3460}.

Ph.\ Barbe, W.P.\ McCormick (2010). Invariance principle for some FARIMA and
nonstationary linear processes in the domain of a stable law, 
{\tt arxiv:1007.0576}.

B.\ Bassan, E.\ Bona (1990). Moments of stochastic processes governed by
Poisson random measures, {\sl Comment.\ Math.\ Univ.\ Carolinae}, 31, 337--343.




N.H.\ Bingham, C.M.\ Goldie, J.L.\ Teugels (1989). {\sl Regular Variation},
2nd ed. Cambridge University Press.





M.\ Cs\"org\H o, S.\ Cs\"org\H o, L.\ Horv\`ath, D.\ Mason (1986). Normal and
stable convergence of integral functions of the empirical distribution
function, {\sl Ann.\ Probab.}, 14, 86--118.



M.\ Donsker (1951). An invariance principle for certain probability limit
theorem, {\sl Four papers on probability}, {\sl Mem.\ Amer.\ Math.\ Soc.}, 6.


J.\ Geffroy (1958/1959). Contribution \`a la th\'eorie des valeurs extr\^emes,
{\sl Publications de l'Institut de Statistique des Universit\'es de Paris},
7--8, 37--185.

S.\ Ghosh, G.\ Samorodnitsky (2009). The effect of memory on functional large
deviations of infinite moving average processes, {\sl Stoch.\ Proc. Appl.},
119, 534--561.

C.W.J.\ Granger (1980). Long memory relationships and the aggregation of
dynamic models, {\sl J.\ Econometrics}, 14, 227--238.

C.W.J.\ Granger (1988). Models that generate trends, {\sl J.\ Time Series
Anal.}, 9, 329--343.




J.\ Kiefer (1972). Iterated logarithm analogues for sample quantiles when 
$p_n\downarrow 0$, {\sl Proc.\ Sixth Berkeley Sympos.\ on Math.\ Statist.\ and
Probab.}, 1, 227--244.

J.F.C.\ Kingman (1961). The single server queue in heavy traffic, {\sl Proc.\
Camb.\ Phil.\ Soc.}, 57, 902--904.

J.F.C.\ Kingman (1962). On queues in heavy traffic, {\sl J.\ Roy.\ Statist.\
Soc., ser.\ B}, 24, 383--392.

J.F.C.\ Kingman (1965). The heavy traffic approximation in the theory of 
queues, in {\sl Proc.\ Symp.\ on Congestion Theory}, W.L.\ Smith and W.E.\ 
Wilkinson eds., pp.137--159, University of North Carolina Press.

K.M.\ Kosi\'nski, O.J.\ Boxma, B.\ Zwart (2010). Convergence of the all-time
supremum of a L\'evy process in the heavy-traffic regime, 
{\tt arXiv:1007.0155v1}.



P.C.B.\ Philipps (1987). Time series with a unit root, {\sl Econometrica},
55, 277--301.







N.\ Privault (2009). Moment identities for Poisson-Skorokhod integrals and
application to measure invariance, {\it C.\ R.\ Acad.\ Sci.\ Paris, S\'er.\ I},
347, 1071--1074.

N.\ Privault (2010). Invariance of Poisson measures under random 
transformations, {\tt arXiv:1004.2588v1}.

Yu.\ V.\ Prohorov (1963). Transition phenomena in queueing processes, I. 
(Russian), {\sl Litovsk.\ Mat.\ Sb.}, 3, 199--205.

Ph.\ Protter (1992). {\sl Stochastic Integration and Differential Equations,
a New Approach}, Springer.

W.\ Pruitt (1981). The growth of random walks and L\'evy processes, {\sl 
Ann.\ Probab.}, 9, 948--956.



S.I.\ Resnick (2007). {\sl Heavy-Tail Phenomena, Probabilistic and Statistical
Modeling}, Springer.


S.\ Shneer, V.\ Wachtel (2009). Heavy-traffic analysis of the maximum of
an asymptotically stable random walk, {\tt arXiv:0902.2185}.

G.R.\ Shorack, J.A.\ Wellner (1978). Linear bounds on the empirical 
distribution function, {\sl Ann.\ Probab.}, 6, 349--353.

G.R.\ Shorack, J.A.\ Wellner (1986). {\sl Empirical Processes with Applications
to Statistics}, Wiley.



W.\ Szczotka, W.A.\ Woykzy\'nski (2003). Distribution of suprema of L\'evy 
processes via the heavy traffic invariance principle, {\sl Prob.\ Math.\ 
Statist.}, 23, 251-272.

S.R.S.\ Varadhan (1966). Asymptotic probabilities and differential equations,
{\sl Comm.\ Pure Appl.\ Math.}, 19, 261--286.

N.\ Veraverbeke (1977). Asymptotic behavior of
Wiener-Hopf factors of a random walk. {\sl Stoch. Proc. Appl.}, 5, 27--37.

\hfuzz=1pt
W.\ Whitt (2002){\sl Stochastic-Process Limits, an Introduction to 
Sto\-chastic-Process Limits and their Applications to Queues}, Springer.

\hfuzz=0pt

W.B.\ Wu, X.\ Shao (2006). Invariance principle for fractionally integrated
nonlinear processes, in {\sl Recent Developments in Nonparmetric Inference
and Probability}, {\sl IMS Lecture Notes-Monographe Series}, 50, 20-30.

}

\bigskip

\setbox1=\vbox{\halign{#\hfil&\hskip 40pt #\hfill\cr
  Ph.\ Barbe            & W.P.\ McCormick\cr
  90 rue de Vaugirard   & Dept.\ of Statistics \cr
  75006 PARIS           & University of Georgia \cr
  FRANCE                & Athens, GA 30602 \cr
  philippe.barbe@math.cnrs.fr                        & USA \cr
                        & bill@stat.uga.edu \cr}}%
\box1%

\bye
\vfill\eject

\noindent {\bf Proof that \hypMg\ implies that \poorBold{$(g_n)$} is normalized
regularly varying.}

\bigskip

Set $h(x)=1/g_{\lfloor x\rfloor}$. Setting $\lambda=i/n$, \hypMg\ yields
$$
  \sup_{n^{-\delta}\leq\lambda\leq 1} \Bigl|{h(n)\over h(\lambda n)}-
  \Bigl({\lfloor \lambda n\rfloor\over n}\Bigr)^{\gamma-1}\Bigr|
  = o(n^{-\delta}) \, .
  \eqno{\equa{gAa}}
$$
Note that in the range $n^{-\delta}\leq\lambda\leq 1$,
$$
  \Bigl| \Bigl({\lfloor n\lambda\rfloor\over n}\Bigr)^{\gamma-1}
  -\lambda^{\gamma-1}\Bigr| 
  \leq c \Bigl|{\lfloor n\lambda\rfloor\over n}-\lambda
    \Bigr|^{(\gamma-1)\wedge 1}
  \leq {c\over n^{(\gamma-1)\wedge 1}} \, .
$$
So, if we choose $\delta$ to be less than $(\gamma-1)\wedge 1$, we have
$$
  \sup_{n^{-\delta}\leq\lambda\leq 1}
  \Bigl| {h(n)\over h(\lambda n)}-\lambda^{\gamma-1}\Bigr| = o(n^{-\delta})
  \, . 
  \eqno{\equa{gA}}
$$
Set $i=\lambda n$ so that $n=i/\lambda$ and $n^{-\delta}\leq \lambda\leq 1$. 
We have 
$$
  n^{-\delta}
  =(i/\lambda)^{-\delta}
  = i^{-\delta} \lambda^\delta
  \leq i^{-\delta} \, .
$$
So \gA\ yields
$$
  \sup_{i^{-\delta}\leq\lambda\leq 1} \Bigl|{h(i/\lambda)\over h(i)}
    -\lambda^{\gamma-1}\Bigr| = o(i^{-\delta}) \, .
$$
Thus, in this range of $\lambda$, and substituting $n$ for $i$,
$$
  {h(n/\lambda)\over h(n)} 
  = \lambda^{\gamma-1} +o(n^{-\delta}) \, . 
$$
Substituting $1/\lambda$ for $\lambda$, we obtain
$$
  {h(\lambda n)\over h(n)} 
  = \Bigl({1\over\lambda}\Bigr)^{\gamma-1}+ o(n^{-\delta})
  =\lambda^{1-\gamma}\bigl( 1+o(n^{-\delta})\lambda^{\gamma-1}\bigr)
$$
in the range $1\leq\lambda\leq n^\delta$.

Next, \gA\ implies
$$
  {h(n)\over h(\lambda n)} 
  = \lambda^{\gamma-1} +o(n^{-\delta})
  = \lambda^{\gamma-1}\bigl( 1+o(n^{-\delta})\lambda^{1-\gamma}\bigr) \, .
  \eqno{\equa{gB}}
$$
We restrict $\lambda$ to a range such that $\lambda^{2(1-\gamma)}n^{-\delta}
\leq 1$, that is, given \gAa, 
$\lambda\geq n^{-\delta/2(\gamma-1)}\vee n^{-\delta}$. Then \gB\ implies
$$\eqalign{
  {h(\lambda n)\over h(n)}
  &{}=\lambda^{1-\gamma}\bigl( 1+o(n^{-\delta})\lambda^{1-\gamma}\bigr)\cr
  &{}=\lambda^{1-\gamma}+o(n^{-\delta/2}) \, . \cr
  }
$$
So, we have
$$
  \sup_{n^{-\delta/2(\gamma-1)}\vee n^{-\delta}\leq\lambda\leq n^\delta}
  \Bigl| {h(\lambda n)\over h(n)} -\lambda^{1-\gamma}\Bigr| = o(n^{-\delta/2})
  \, .
$$
Setting 
$$
  \kappa={\delta\over 2}\wedge {\delta\over 2(\gamma-1)} \, ,
$$
this implies
$$
  \sup_{n^{-\kappa}\leq\lambda\leq n^\kappa}
  \Bigl| {h(n\lambda)\over h(n)}-\lambda^{1-\alpha}\Bigr|
  = o(n^{-\kappa}) \, .
$$
Substituting $h$ for $F^\leftarrow(1-1\Id)$ in Proposition \hypRVMeaning.i,
the second statement of that proposition implies that $h(n)=cn^{1-\gamma}
\bigl(1+o(n^{-\epsilon}\bigr)$ for some positive $\epsilon$.

\bigskip

{\it So I propose to add in the paper that, considering $(1/g_n)$, 
Proposition \hypRVMeaning\ implies that for some positive $\epsilon$
$$
  {1\over g_n}=c n^{1-\gamma}\bigl(1+o(n^{-\epsilon})\bigr)
$$
and so $(g_n)$ is normalized regularly varying.}

\bye

mpage -2 -c -o -M-100rl-80b-220t -t veraverbeke7.ps > toto.ps

\bye